\documentclass[english]{article}

\usepackage{geometry}
\usepackage[T1]{fontenc}
\usepackage[utf8]{inputenc}
\usepackage{bm}
\usepackage{graphicx}
\graphicspath{{../imgs/}}
\usepackage{tikz}
\usetikzlibrary{positioning, arrows.meta, calc}
\tikzset{
  box/.style={
    draw,
    rounded corners,
    minimum width=2cm,
    minimum height=1cm,
    align=center,
  },
  arrow/.style={
    ->,
    >=Latex,
    thick,
    draw=blue,
  },
  dashed-line/.style={
    dashed,
    thick,
    draw=blue,
  }
}

\usepackage{caption}
\usepackage{subcaption}
\usepackage{amssymb}
\usepackage{mathrsfs}
\usepackage{amsmath}
\usepackage{dsfont}
\usepackage{extarrows}
\usepackage{enumitem}
\usepackage{booktabs}
\usepackage{makecell}
\usepackage{tablefootnote}
\usepackage{threeparttable}
\usepackage{tabularx}
\usepackage{multirow}
\usepackage[linesnumbered,ruled,vlined]{algorithm2e}
\usepackage{xcolor}

\SetCommentSty{mycommfont}

\SetKwInput{KwInput}{Input}                
\SetKwInput{KwOutput}{Output}              

\usepackage[unicode=true,
 bookmarks=false,
 breaklinks=false,pdfborder={0 0 1},colorlinks=false]
 {hyperref}
\hypersetup{
 colorlinks,citecolor=blue,filecolor=blue,linkcolor=blue,urlcolor=blue}

\makeatletter
\usepackage{amsthm}
\usepackage{comment}
\usepackage{natbib}

\graphicspath{ {./imgs/} }

\newtheorem{theorem}{Theorem}[section]
\newtheorem{proposition}[theorem]{Proposition}

\newtheorem{assumption}[theorem]{Assumption}

\theoremstyle{remark}
\newtheorem{remark}[theorem]{Remark}
\newcounter{theoremctr}
\newcounter{corollaryctr}
\newcounter{assumptionctr}
\newcounter{propositionctr}
\newcounter{lemmactr}

\renewcommand{\thetheoremctr}{\arabic{theoremctr}}     
\renewcommand{\thecorollaryctr}{\arabic{corollaryctr}}  
\renewcommand{\theassumptionctr}{\arabic{assumptionctr}} 
\renewcommand{\thepropositionctr}{\arabic{propositionctr}} 
\renewcommand{\thelemmactr}{\arabic{lemmactr}}

\allowdisplaybreaks

\let\bar\overline
\let\tilde\widetilde

\title{Berry-Esseen bounds for multivariate martingale difference sequences in the Kolmogorov distance}

\author{
\begin{tabular}{c}
Weichen Wu\thanks{The Voleon Group, New York, NY 10010, USA.}
\quad
Dung Le\thanks{Department of Statistics and Data Sciences, University of Texas, Austin, TX 78705, USA.}
\quad
Arun Kumar Kuchibhotla\thanks{Department of Statistics and Data Science, Carnegie Mellon University, Pittsburgh, PA 15213, USA.}
\quad
Alessandro Rinaldo\footnotemark[2]
\end{tabular}
}

\date{\today}

\begin{document}

\maketitle

\begin{abstract}
    We derive new Gaussian approximations for finite martingale difference sequences in $\mathbb{R}^d$ with respect to the Kolmogorov distance. Under appropriate conditions, our bounds exhibit a dependence of order $n^{-1/4}$  on the length of the sequence and of order $\mathrm{polylog}(d)$ on the dimension. As an application, we derive a high-dimensional Berry-Esseen bound over hyper-rectangles for martingale sequences generated from Markov chains.
\end{abstract}

\section{Introduction}

In this note, we establish novel high-dimensional Berry-Esseen bounds for $d$-dimensional martingale difference sequences with deterministic terminal quadratic variation over the class of $d$-dimensional hyper-rectangles under various moment conditions. Concretely, let $\{\bm{X}_k\}_{k=0,\ldots,n}$ be a martingale difference sequence in $\mathbb{R}^d$ (thus, $\bm{X}_0 = 0$) having bounded third moments with respect to the filtration $\{\mathscr{F}_k\}_{k=0,\ldots,n}$, with $\mathscr{F}_0$ the trivial $\sigma$-field\footnote{Our results can be easily extended to the case in which $\mathscr{F}_0$ is not trivial.}. Throughout, we will assume that the terminal quadratic variation of the martingale sequence is deterministic. 
\begin{assumption}\label{as:fixed-variation}
It holds that, almost surely,
\begin{equation}\label{eq:uno}
\sum_{k=1}^n \mathbb{E}[\bm{X}_k\bm{X}_k^\top \mid \mathscr{F}_{k-1}] = \sum_{k=1}^n \mathbb{E}[\bm{X}_k\bm{X}_k^\top \mid \mathscr{F}_0]  = n\bm{\Sigma}_n, 
\end{equation}
for a positive definite matrix $\bm{\Sigma}_n$. 
\end{assumption}

As explained below, though seemingly strong, this assumption is widely used in the literature on Gaussian approximation for martingale sequences.


We seek to establish a bound on Kolmogorov distance between the normalized sum of the martingale difference terms (thus a centered martingale) and a centered Gaussian vector with covariance $\bm{\Sigma}_n$. 
That is, we seek to bound
\begin{equation}
    \label{eq:due}
    d_{\mathsf{K}}\left(\frac{1}{\sqrt{n}}\sum_{k=1}^n \bm{X}_k, \frac{1}{\sqrt{n}}\bm{T}_n \right) : = \sup_{\mathcal{A} \in \mathscr{H}_d} \left|\mathbb{P}\left(\frac{1}{\sqrt{n}}\sum_{k=1}^n \bm{X}_k \in \mathcal{A} \right) - \mathbb{P}\left( \frac{1}{\sqrt{n}}\bm{T}_n \in \mathcal{A} \right)\right|,
\end{equation}
where $\mathscr{H}_d = \{ \prod_{i=1}^d [a_i, b_i], -\infty <  a_i \leq b_i < \infty \}$ denotes the set of finite hyper-rectangles in $\mathbb{R}^d$ and $\frac{1}{\sqrt{n}} \bm{T}_n \sim \mathcal{N}(\bm{0},\bm{\Sigma}_n)$.


The literature on Berry--Esseen bounds for univariate martingale difference sequences is extensive, with many of the key results established under the deterministic terminal quadratic variation assumption \eqref{eq:uno} or stronger versions of it. We refer the reader to \cite{hall2014martingale}, \cite{grams1972rates}, \cite{kato1979rates}, \cite{bothausen1982convergence}, \cite{renz1996note}, \cite{rinott2000bounds}, \cite{el2007exact},  \cite{mourrat2013rate}, \cite{Fan2019Exact}, \cite{wu2020berry} and references therein, for various results, assumptions and discussion. Specifically, as originally shown in \cite{grams1972rates} and \cite{bothausen1982convergence}, for martingale difference sequences with uniformly bounded third moments and deterministic conditional variances -- a more stringent condition than \eqref{eq:uno} --  the optimal rate for \eqref{eq:due} is of order $O\left( \frac{\log n }{n^{1/4}} \right)$. Faster rates of order $O\left( \frac{\log n }{n^{1/2}} \right)$ are possible with additional conditions, such as under boundedness and  higher order moment assumptions: see, e.g., \cite{bothausen1982convergence}, \cite{Haausler:88}, \cite{renz1996note}, \cite{rinott2000bounds},  \cite{Fan2019Exact} and \cite{wu2020berry}. In general, as noted by  \cite{rinott2000bounds}, assuming only a deterministic terminal quadratic variation as in \eqref{eq:uno} and existence of third moments, faster rates than $O(n^{-1/8})$ are not possible.

As indicated above, the deterministic quadratic variation condition \eqref{eq:uno} is commonly found in the literature on univariate Berry--Esseen bounds for martingales as well as in the more recent contributions on high-dimensional Gaussian approximations (see, e.g., \cite{JMLR2019CLT} and \cite{wu2025uncertainty}). 
While mathematically very convenient, Assumption \ref{as:fixed-variation} is, of course, quite restrictive and unrealistic for most settings.
It can be relaxed at the expense of an additional term in the Berry--Esseen bound that depends on the magnitude of the difference between $\sum_{k=1}^n \bm{V}_k $ and $n\bm{\Sigma}_n$. Indeed, it is known (see, e.g., \cite{hall2014martingale}) that for martingale difference sequences, the accuracy of the Gaussian approximation is typically driven by this term. As this type of generalization relies on well-known arguments, in order to avoid technicalities and focus only on the error of the Gaussian approximation, we will present our results assuming that \eqref{eq:uno} holds. In Section~\ref{sec:application} 
we outline a standard strategy (see, e.g., \cite{Dvoretzky.clt.martingale:72,belloni2018highdimensionalcentrallimit} and \cite{yurinksi.martingale}), to extend our results to martingales with non-deterministic terminal quadratic variation and exemplify this approach by deriving a novel Berry--Esseen bound over hyper-rectangles for a special class of martingale difference sequences arising from reinforcement learning and the study of Markov chains.

In contrast to the wealth of results available in the univariate case, the literature on high-dimensional Gaussian approximation of martingale difference sequences is comparatively limited. The multivariate problem is substantially harder. Beyond the added complexity of working in $\mathbb{R}^d$, the choice of probability metric matters: the Wasserstein and convex distances admit comparatively clean analyses, but the Kolmogorov distance — defined as the supremum over all hyper-rectangles — is both the most difficult to control and the most directly relevant to statistical applications such as simultaneous inference and the construction of joint confidence regions.
 \cite{belloni2018highdimensionalcentrallimit} derived a bound for smooth (thrice differentiable) functions of martingale difference sequences, while \cite{JMLR2019CLT} considered instead functions that are only twice differentiable, assuming the terminal quadratic variation condition \eqref{eq:uno}. By extending some arguments in \cite{belloni2018highdimensionalcentrallimit}, \cite{yurinksi.martingale} derived general bounds for high-dimensional approximations of martingale difference sequences over arbitrary classes of subsets of $\mathbb{R}^d$, including hyper-rectangles. Despite the generality of this result, the resulting rates are sub-optimal in $n$.  More recently, \cite{srikant2024rates} and 
\cite{wu2025uncertainty} obtained high-dimensional Gaussian approximation rates in the Wasserstein distance (the latter also under Assumption~\ref{as:fixed-variation}), from which it is possible to also derive bounds in the convex distance (see Proposition A.1 in \cite{nourdin2021multivariate}).

To the best of our knowledge there are only two recent contributions in the literature, offering high-dimensional Berry-Esseen martingale bounds over hyper-rectangles.
\cite{Kojevnikov2022}, building on the proof technique by \cite{arun2020hdCLT}, derived a high-dimensional Berry—Esseen bound  in the Kolmogorov distance \eqref{eq:due} of order $O \left( (\log d)^{5/4}  \frac{\log n}{n^{1/4}} \right)$ for martingale difference sequences in $\mathbb{R}^d$. However, their result requires that each conditional covariance matrix $\mathbb{E}[\mathbf{X}_k \mathbf{X}_k^\top \mid \mathcal{F}_{k-1}]$
 be measurable with respect to the initial $\sigma$-field $\mathcal{F}_0$, a very stringent condition, stronger than the deterministic terminal quadratic condition \eqref{eq:uno}, that is not satisfied in most examples arising from Markov chains or stochastic approximation algorithms. Very recently, \cite{rubtsov2026gaussian}, in their Theorem 4, gave a novel high-dimensional Gaussian approximation bound for martingale difference sequences under Assumption~\ref{as:fixed-variation}. The authors' approach combines a generalization of the proof technique used in Theorem 2.1 of \cite{rollin2018}, establishing a Berry-Esseen theorem in the Wasserstein distance for univariate martingale difference sequences using Stein's method, with recent sharpening of the properties of the solution to the multivariate Stein equation, as in \cite{gallouet2018regularity} and \cite{wu2025uncertainty}. Their bound holds under minimal assumptions and depends on the choice of a $d$-dimensional positive definite matrix $\Sigma$, acting as a spectral regularizer. Under additional assumptions, it leads to high-dimensional martingale Berry-Esseen bounds without requiring a deterministic terminal quadratic variation. We will compare our results with those of \cite{Kojevnikov2022} and \cite{rubtsov2026gaussian} below. 




 In this paper, we establish some new multivariate Berry–Esseen bounds in Kolmogorov distance assuming the fixed terminal quadratic variation condition \eqref{eq:uno} as well as appropriate moment conditions. 
 Under a moment ratio condition (see \eqref{eq:telescope-condition} below), in Theorem~\ref{thm:Berry-Esseen-Kol} we obtain an $O(n^{-1/4} \log n)$ bound for the Kolmogorov distance \eqref{eq:due}, and with direct moment value conditions (see \eqref{eq:MVB-upper}) we achieve the optimal $O(n^{-1/4})$ rate in Theorem \ref{thm:Berry-Esseen-Kol-MVB}. Our bounds match the form of \cite{Kojevnikov2022} while imposing a strictly weaker structural assumption. We further demonstrate that the $n^{-1/4}$ rate is optimal under finite third moments by a lower bound construction based on \cite{bothausen1982convergence}. Without explicit moment bound conditions (though still assuming finite third moments) we are only able to guarantee a rate of order $n^{-1/8}$; see Theorem \ref{thm:Berry-Esseen-Kol-naive}.  All our bounds exhibit only a polylogarithmic dependence on the dimension $d$, as expected.  As an application, and dropping the deterministic terminal quadratic variation assumption \eqref{eq:uno}, we derive Kolmogorov-distance Berry–Esseen bounds for additive functionals of Markov chains, extending the Wasserstein-distance results of \cite{wu2025uncertainty} to the stronger distance.


{\bf Notation.} Throughout the paper, for a square matrix $\bm{A}$, we use $\lambda_{\max}(\bm{A}), \lambda_{\min}(\bm{A})$ to denote its maximum and minimum eigenvalues; $d_{\max}(\bm{A}), d_{\min}(\bm{A})$ to denote its maximum and minimum diagonal elements; and $\|\bm{A}\|$ its operator norm. If $\bm{A}$ is positive semi-definite, then $\|\bm{A}\| = \lambda_{\max}(\bm{A})$.
 We will write $\log_+d := \max\{1,\log d\}$ and, for any positive integer $n$, $[n] = \{ 1,2,\ldots,n\}$.

\section{Problem setting and Results}\label{sec:setting}
Let $\{\bm{X}_k\}_{k =0,1,\ldots,}$, where $\bm{X}_0 = 0$, be a martingale difference process in $\mathbb{R}^d$ with respect to the $\sigma-$field filtration $\{\mathscr{F}_k\}_{0 \leq k \leq n}$, with $\mathcal{F}_0$  the trivial $\sigma$-field. Set, for each $k=1,\ldots,n$, $\bm{S}_k = \sum_{i=1}^k \bm{X}_i$, so that, $\{\bm{S}_k\}_{k =1,\ldots,n}$ is a $d$-dimensional centered martingale. We are interested in deriving a high-dimensional Gaussian approximation to the normalized sum $\bm{S}_n/\sqrt{n}$, under the standard Assumption \ref{as:fixed-variation}. We derive such an approximation by bounding the $d$-dimensional Kolmogorov distance between $\bm{S}_n/\sqrt{n}$ and the normalized sum of centered Gaussian vectors with matching conditional covariances. In detail, let 
\begin{align}\label{eq:defn-Sigma-n}
\bm{\Sigma}_n = \frac{1}{n} \sum_{k=1}^n \mathbb{E}[\bm{X}_k\bm{X}_k^\top \mid \mathscr{F}_0], \quad  \quad \bm{V}_k = \mathbb{E}[\bm{X}_k\bm{X}_k^\top \mid \mathscr{F}_{k-1}], \quad k=1,\ldots,n,
\end{align}
and $ \frac{1}{\sqrt{n}} \bm{T}_n : = \frac{1}{\sqrt{n}}\sum_{k=1}^n \bm{V}_k^{1/2} \bm{Z}_k$, where  $\{\bm{Z}_k\}_{k=1,\ldots,n}$ is a sequence of i.i.d. $d$-dimensional standard Gaussian random vectors.
We will derive novel bounds for the Kolmogorov distance \eqref{eq:due} in the introduction, which we can rewrite more explicitly as
\begin{align*}
d_{\mathsf{K}}\left(\frac{1}{\sqrt{n}} \bm{S}_n,\frac{1}{\sqrt{n}}\bm{T}_n \right)  & = \sup_{\mathcal{A} \in \mathscr{H}_d} \left|\mathbb{P}\left(\frac{1}{\sqrt{n}}\sum_{k=1}^n \bm{X}_k \in \mathcal{A} \right) - \mathbb{P}\left(\frac{1}{\sqrt{n}}\sum_{k=1}^n \bm{V}_k^{1/2} \bm{Z}_k \in \mathcal{A} \right)\right|\\
& : = d_{\mathsf{K}}\left(\frac{1}{\sqrt{n}}\sum_{k=1}^n \bm{X}_k,\mathcal{N}(\bm{0},\bm{\Sigma}_n)\right),
\end{align*}
where in the last expression we use the fact that $\frac{1}{\sqrt{n}} \bm{T}_n \sim \mathcal{N}(\bm{0},\bm{\Sigma}_n)$, which follows from the deterministic terminal quadratic variation assumption \eqref{eq:uno}. The first identity above elucidates the scope of the Gaussian approximation, whereby the original martingale difference sequence $\{ \bm{X}_k \}_{k=0,\ldots,n}$ is compared with a new martingale difference sequence $\{ \bm{V}_k^{1/2} \bm{Z}_k\}_{k=0,\ldots,n}$ with matching conditional second moments and conditionally Gaussian; i.e., for every $k$,
\[
\mathrm{Var}[\bm{X}_k \mid \mathscr{F}_{k-1}] = \mathrm{Var}[\bm{V}_k^{1/2} \bm{Z}_k \mid \mathscr{F}_{k-1}] \quad \text{and} \quad \bm{V}_k^{1/2} \bm{Z}_k \mid \mathscr{F}_{k-1} \sim \mathcal{N}(\bm{0},\bm{V}_k).
\]

In our main results below, we present three bounds for $d_{\mathsf{K}}\left(\frac{1}{\sqrt{n}}\sum_{k=1}^n \bm{X}_k,\mathcal{N}(\bm{0},\bm{\Sigma}_n)\right)$ under Assumption~\ref{as:fixed-variation} and  various third moment conditions in the $\| \cdot \|_{\infty}$ norm.

In our proofs, we use various techniques and results from the literature on Gaussian approximations: Lindeberg's swapping, smoothing via the Ornstein-Uhlenbeck process or Gaussian noise, and the explicit form of the solution to Stein's multivariate equation in terms of the generator of the Ornstein-Uhlenbeck semi-group. The effectiveness of these tools in deriving high-dimensional Berry-Esseen bounds in the Kolmogorov distance was demonstrated in \cite{fang2021highdimensional} and \cite{arun2020hdCLT}, and our proofs build on these techniques. For more details, we refer the reader to, e.g.,  \cite{Barbour1990SteinsMF}, \cite{gotze1991rate},   \cite{Bhattacharya2010AnEO}, \cite{meckes2009stein}, and \cite{nourdin.peccati.2012}. 




\subsection{High-dimensional Berry-Esseen bounds in Kolmogorov distance for martingales with fixed terminal quadratic variation}


In our first general bound, in addition to the terminal quadratic variation condition in Assumption~\ref{as:fixed-variation}, we impose the mildest possible third moment condition, resulting in a slow rate in $n$ of order $O(n^{-1/8})$. This is consistent with analogous findings in the univariate literature: see, e.g., \cite{rinott2000bounds}.
\begin{assumption}\label{as:finite-3rd-moment}
For all $k \in [n]$, $\mathbb{E}\|\bm{X}_k\|_{\infty}^3 < \infty$.
\end{assumption}

The following theorem gives a finite-sample Gaussian approximation under these minimal assumptions.

\begin{theorem}\label{thm:Berry-Esseen-Kol-naive}
Let $\{\bm{X}_k\}_{k \in [n]}$ denote a martingale difference process in $\mathbb{R}^d$ that satisfies Assumptions \ref{as:fixed-variation} and \ref{as:finite-3rd-moment}. Then
\begin{align*}
d_{\mathsf{K}}\left(\frac{1}{\sqrt{n}}\sum_{k=1}^n \bm{X}_k,\mathcal{N}(\bm{0},\bm{\Sigma}_n)\right) \lesssim  \left(\frac{1}{\lambda_{\min}(\bm{\Sigma}_n)}\right)^{\frac{3}{8}}\left( \sum_{k=1}^n \frac{\mathbb{E}\|\bm{X}_k\|_{\infty}^3}{n}\right)^{1/4} \left(\frac{d_{\max}(\bm{\Sigma}_n)}{d_{\min}(\bm{\Sigma}_n)}\right)^{\frac{3}{8}}(\log_+ d)^{9/8} n^{-\frac{1}{8}}.
\end{align*}
\end{theorem}


If we further impose a \emph{ratio bound condition} on the third moment, we can improve the dependence on $n$ to $O(n^{-\frac{1}{4}}\log n)$, which, as shown below, is nearly optimal under the assumed conditions. The proof is given in Appendix \ref{app:proof-thm-Berry-Esseen-Kol}. Conditions of this type have been considered in the literature; see, e.g., \cite{rollin2018} and \cite{JMLR2019CLT}.


\begin{theorem}\label{thm:Berry-Esseen-Kol}
Let $\{\bm{X}_k\}_{k \in [n]}$ denote a martingale difference process in $\mathbb{R}^d$ that satisfies Assumption \ref{as:fixed-variation}, and let $\bm{V}_k, \bm{\Sigma}_n$ be defined as in \eqref{eq:defn-Sigma-n}. Suppose that there exists a positive constant $M$ such that for all $k \in [n]$,
\begin{align}\label{eq:telescope-condition}
\mathbb{E}[\|\bm{X}_k\|_{\infty}^3\mid\mathscr{F}_{k-1}] \leq M \lambda_{\min}(\bm{V}_k) \quad \text{almost surely.}
\end{align}
Then
\begin{align}\label{eq:Berry-Esseen-Kol}
d_{\mathsf{K}}\left(\frac{1}{\sqrt{n}}\sum_{k=1}^n \bm{X}_k,\mathcal{N}(\bm{0},\bm{\Sigma}_n)\right) \lesssim \Delta \log_+d \log \frac{1+\Delta^2}{\Delta^2}
\end{align}
where 
\begin{align}\label{eq:Delta}
\Delta = \sqrt{M}\left(\frac{1}{\lambda_{\min}(\bm{\Sigma}_n)}\right)^{\frac{1}{4}} \left(\frac{d_{\max}(\bm{\Sigma}_n)}{d_{\min}(\bm{\Sigma}_n)}\right)^{\frac{1}{4}}\log_+^{\frac{3}{4}} d \cdot n^{-\frac{1}{4}}.
 \end{align}
\end{theorem}

\begin{remark}
In the one-dimensional setting, the moment ratio condition \eqref{eq:telescope-condition} reduces to
\begin{align*}
\mathbb{E} \left[|X_k|^3 \mid \mathscr{F}_{k-1} \right] \leq M \mathbb{E}[X_k^2\mid\mathscr{F}_{k-1}] \quad \text{almost surely, for all }k \in [n].
\end{align*}
Notice that this is equivalent to the assumption outlined in Theorem 1 of \cite{el2007exact} and weaker than the one made by Corollary 2.3 in \cite{rollin2018}, which also assumes a universal bound on the conditional third moment of $|X_k|$. For Berry--Esseen bounds measured by the Wasserstein distance and the convex distance, \cite{wu2025uncertainty} made the analogous assumption that, for every unit vector $\bm{e} \in \mathbb{R}^d$ and every $k \in [n]$,
\begin{align}\label{eq:MRB-convex}
\mathbb{E}[(\bm{e}^\top \bm{X}_k)^2\|\bm{X}_k\|_2\mid \mathscr{F}_{k-1}] \leq M\mathbb{E}[(\bm{e}^\top \bm{X}_k)^2\mid \mathscr{F}_{k-1}] \quad \text{almost surely}.
\end{align}
Notably, the condition \eqref{eq:MRB-convex} can be implied by the assumption that $\|\bm{X}_k\|_2 \leq M$ almost surely; but our condition \eqref{eq:telescope-condition} \emph{cannot} be implied by the assumption that $\|\bm{X}_k\|_{\infty} \leq M$. 
\end{remark}


If instead of the moment ratio bound more direct third \emph{moment value bounds} on the terms of the martingale difference sequence are imposed, then the following theorem can further tighten the bound to $O(n^{-\frac{1}{4}})$. See  Appendix \ref{app:proof-thm-Berry-Esseen-Kol-MVB} for a proof.

\begin{theorem}\label{thm:Berry-Esseen-Kol-MVB}
Let $\{\bm{X}_k\}_{k \in [n]}$ denote a martingale difference sequence in $\mathbb{R}^d$ that satisfies Assumption \ref{as:fixed-variation}, and let $\bm{V}_k, \bm{\Sigma}_n$ be defined as in \eqref{eq:defn-Sigma-n}. Suppose that there exist constants $\beta>0,\gamma>0$, such that for all $k \in [n]$,
\begin{align}\label{eq:MVB-upper}
\mathbb{E}\|\bm{X}_k\|_{\infty}^3 \leq (\gamma \log_+ d)^{3/2}, \quad \text{and} \quad d_{\max}(\bm{V}_k) \leq \beta \quad \text{almost surely}.
\end{align}
It can then be guaranteed that
\begin{align}\label{eq:Berry-Esseen-Kol-MVB1}
d_{\mathsf{K}}\left(\frac{1}{\sqrt{n}}\sum_{k=1}^n \bm{X}_k,\mathcal{N}(\bm{0},\bm{\Sigma}_n)\right) \lesssim  \left(\frac{\gamma + \beta}{d_{\min}(\bm{\Sigma}_n)}\right)^{\frac{3}{8}} \left(\log_+^{\frac{3}{2}}d\right) n^{-\frac{1}{8}}.
\end{align}
If we further assume that there exists $\alpha>0$ such that 
\begin{align}\label{eq:MVB-lower}
\min_{k=1,\ldots,n} \lambda_{\min}(\bm{V}_k) \geq \alpha, \quad \text{almost surely},
\end{align}
then 
\begin{align}\label{eq:Berry-Esseen-Kol-MVB2}
d_{\mathsf{K}}\left(\frac{1}{\sqrt{n}}\sum_{k=1}^n \bm{X}_k,\mathcal{N}(\bm{0},\bm{\Sigma}_n)\right) \lesssim  \frac{(\gamma + \beta)^{3/4}}{\sqrt{\alpha d_{\min}^{1/2}(\bm{\Sigma}_n)}} \left(\log_+^{2}d\right) n^{-\frac{1}{4}}.
\end{align}
\end{theorem}


Several remarks are in order.
\begin{itemize}
\item  The upper bound \eqref{eq:Berry-Esseen-Kol-MVB2} takes exactly the same form as the one given in the Remark following Theorem 3.1 of \cite{Kojevnikov2022}. The main notable difference is that our result does not depend on every conditional covariance matrix $\bm{V}_k$ to be measurable with respect to $\mathscr{F}_0$, i.e. to be deterministic, for all $k$; in fact, we only require $\bm{\Sigma}_n$ to be fixed.
A second minor difference is that \cite{Kojevnikov2022} achieved a better dependence on $d$, of order $(\log_+ d)^{\frac{5}{4}}$ instead of $(\log_+ d)^{2}$.
\item Notice that Theorem \ref{thm:Berry-Esseen-Kol-MVB} directly implies that, for bounded martingale difference sequences with $\max_k \| \bm{X}_k \|_\infty \leq R$ almost surely,
\begin{align*}
d_{\mathsf{K}}\left(\frac{1}{\sqrt{n}}\sum_{k=1}^n \bm{X}_k,\mathcal{N}(\bm{0},\bm{\Sigma}_n)\right) \lesssim  \frac{R^{3/2}}{\sqrt{\alpha d_{\min}^{1/2}(\bm{\Sigma}_n)}} \left(\log_+^{2}d\right) n^{-\frac{1}{4}}.
\end{align*}
\item Right before posting this paper as a preprint, we became aware of the very recent and interesting  work of \cite{rubtsov2026gaussian}, which established that for every $\bm{\Sigma} \succeq 0$,
\begin{align*}
d_{\mathsf{K}}\left(\frac{1}{\sqrt{n}}\sum_{k=1}^n \bm{X}_k,\mathcal{N}(\bm{0},\bm{\Sigma}_n)\right) \lesssim \frac{\log_+^{\frac{5}{4}} d}{\sqrt{d_{\min}^{1/2}(\bm{\Sigma}_n)}}\left(d_{\max}(\bm{\Sigma}) + \sum_{k=1}^n \mathbb{E} \frac{\|\bm{X}_k\|_{\infty}^3}{\lambda_{\min}(\bm{P}_k + \bm{\Sigma})}\right)^{1/2} n^{-\frac{1}{4}}.
\end{align*}
Here, $\bm{P}_k = \sum_{i=k}^n \mathbb{E}[\bm{X}_i\bm{X}_i^\top \mid \mathscr{F}_{i-1}]$. One can verify that under conditions \eqref{eq:MVB-upper} and \eqref{eq:MVB-lower}, this bound translates to a rate scaling as
\begin{align*}
\frac{\gamma^{3/4}}{\sqrt{\alpha d_{\min}^{1/2}(\bm{\Sigma}_n)}} \left(\log_+^{2}d\right) n^{-\frac{1}{4}} \log^{\frac{1}{2}} n.
\end{align*}
\end{itemize}



    Under the assumptions \eqref{eq:MVB-upper} and \eqref{eq:MVB-lower}, the convergence rate $O(n^{-1/4})$ for the Kolmogorov distance achieved by our theorem is not improvable. In particular, based on the construction presented in Section 6 of \cite{bothausen1982convergence}, one can construct an example for which the Kolmogorov distance converges at the rate $O(n^{-1/4})$. We define a triangular array of random vectors  $\bm{X}_{ni}$ ($1\leq i\leq n$) in $\mathbb{R}^{d}$. For $2\leq j \leq d$, the $j$-coordinate random variables $X^j_{ni}$ are i.i.d. $\mathcal{N}(0,1)$. For the random variables in the first coordinate, we follow the same construction as the Markov chain in Example 1 of \cite{bothausen1982convergence}. Concretely, for $1\leq i\leq \lfloor n - 2\sqrt{n}\rfloor$ or $\lfloor n - \sqrt{n}\rfloor < i \leq n$, $X^1_{ni}$ are i.i.d. standard normal random variables. For $\lfloor n - 2\sqrt{n}\rfloor < i\leq \lfloor n - \sqrt{n}\rfloor$, the distribution of $X_{ni}^1$ is defined conditionally on $S_{n,i-1} := \sum_{j=1}^{i-1}X_{nj}^1$: 
    
    \begin{equation*}
       X^1_{ni}|S_{n,i-1} = x\sim \begin{cases}
      \frac{1}{2} (\delta_{-1}+\delta_{1}), & \text{ if  } x \notin [-\sqrt{n}\lambda_i/4,\sqrt{n}\lambda_i/4];\\
    16\lambda_i^2 \delta_{\rho_i} + (1-16\lambda_i^2)\delta_{-1/\rho_i}, &\text{ if } x \in [-\sqrt{n}\lambda_i/4,\sqrt{n}\lambda_i/4].
        \end{cases}
    \end{equation*}
    Here, $\lambda_i = \left(1-\dfrac{i}{n}\right)^{1/2}$, $\rho_i = \left(\dfrac{1-16\lambda_i^2}{16\lambda_i^2} \right)^{1/2}$, and $\delta_x$ denotes the Dirac measure on $x$. Then, thanks to the independence of each component in these random vectors $\bm{X}_k$, it is straightforward to verify that $\{\bm{X}_{ni}\}_{1\leq i\leq n}$ is a martingale difference sequence satisfying 
    \begin{equation*}
        \bm{V}_{nk} = \mathbb{E}[\bm{X}_{nk}\bm{X}_{nk}^\top|\mathscr{F}_{k-1}] = \bm{I}_d = \mathbb{E}[\bm{X}_{nk}\bm{X}_{nk}^\top|\mathscr{F}_{0}]. 
    \end{equation*}
    This means Assumption \ref{as:fixed-variation} is satisfied with $\bm{\Sigma}_n = \bm{I}_d$. In addition, Proposition 1 of \cite{bothausen1982convergence} states that there exists a universal constant $C$ such that for all $1 \leq k \leq n$, $\mathbb{E}|X_{nk}^1|^3 \leq C$. Therefore, the third moment of $\|\bm{X}_{nk}\|_{\infty}$ can be bounded by
    \begin{align*}
    \mathbb{E}\|\bm{X}_{nk}\|_{\infty}^3 &= \mathbb{E}[\max\{|X_{nk}^1|, \|\bm{X}_{nk}^{(-1)}\|_{\infty}\}]^3 \\
    &\leq \mathbb{E}[|X_{nk}^1| + \|\bm{X}_{nk}^{(-1)}\|_{\infty}]^3 \\ 
    &\overset{\text{a}}{\leq} 4 \mathbb{E}|X_{nk}^1|^3 + 4\mathbb{E}\|\bm{X}_{nk}^{(-1)}\|_{\infty}^3 \\ 
    &\overset{\text{b}}{\leq} 4C + 4 \log_+^{\frac{3}{2}}(d-1) < ((4C+4)\log_+d)^{\frac{3}{2}}.
    \end{align*}
    Here, $\bm{X}_{nk}^{(-1)}$ denotes the vector $(\bm{X}_{nk}^{2},\bm{X}_{nk}^3,...,\bm{X}_{nk}^d)^\top$, which has distribution $\mathcal{N}(\bm{0},\bm{I}_{d-1})$; the inequality (a) follows from the inequality $(x+y)^3 \leq 4x^3 + 4y^3$ for all $x,y>0$, and (b) follows from Theorem \ref{thm:Gaussian-inf-norm} in the appendix. Therefore, the conditions \eqref{eq:MVB-upper} and \eqref{eq:MVB-lower} are satisfied with $\gamma = 4C+4, \beta = 1$ and  $\alpha = 1$, respectively.
    Meanwhile, part (c) of Proposition 1 in \cite{bothausen1982convergence} states that
    \begin{align*}
    \lim \sup_{j \to \infty} \sup_{n \geq i \geq j} i^{\frac{1}{4}}d_{\mathsf{K}} (S_{ni}^1/\sqrt{i},\mathcal{N}(0,1)) > 0.
    \end{align*}
    
    Therefore, for $\bm{S}_{ni}: = \sum_{k=1}^i\bm{X}_{nk}$, we can estimate a lower bound for the Kolmogorov distance:
    \begin{align*}
        &\lim \sup_{j \to \infty} \sup_{n \geq i \geq j} i^{\frac{1}{4}}d_{\mathsf{K}}(\bm{S}_{ni}/\sqrt{i}, \mathcal{N}(\bm{0},\bm{I}_d)\mathbb) \\ 
        &= \lim \sup_{j \to \infty} \sup_{n \geq i \geq j} i^{\frac{1}{4}}\sup_{[a,b]}\left|\mathbb{P}[\bm{S}_{ni}/\sqrt{i} \in R_{a,b}] - \mathbb{P}[\mathcal{N}(\bm{0},\bm{I}_d)\in R_{a,b}]\right| \\
        &\overset{\text{(a)}}{=} \lim \sup_{j \to \infty} \sup_{n \geq i \geq j} i^{\frac{1}{4}} \sup_{[a,b]}\left|\mathbb{P}[{S}^{1}_{ni} \in [a_1,b_1]]-\mathbb{P}[\mathcal{N}(0,1) \in [a_1,b_1]]\right| \prod_{k=2}^d\mathbb{P}[\mathcal{N}(0,1)\in[a_k,b_k]]\\
        &\overset{\text{(b)}}{=} \lim \sup_{j \to \infty} \sup_{n \geq i \geq j} i^{\frac{1}{4}}\sup_{[a_1,b_1]}\left|\mathbb{P}[{S}^{1}_{ni} \in [a_1,b_1]]-\mathbb{P}[\mathcal{N}(0,1) \in [a_1,b_1]]\right| >0,
    \end{align*}
    where (a) is from the independence for each component in our random vectors and the normality of the components except the first one, (b) is true by taking $a_i \to -\infty$ and $b_i \to \infty$ for $i\geq 2$, and (c) is due to part (c) of Proposition 1 in \cite{bothausen1982convergence}.

\subsection{Application to martingales generated from Markov chains}\label{sec:application}


The theoretical results in the previous section are all based on the fixed terminal quadratic variation Assumption \ref{as:fixed-variation}, which is usually too strong in practice. In the following theorem, we present a Berry--Esseen bound where the multivariate martingale is generated from a Markov chain. Sequences of this type are widely used in the analysis of MCMC and stochastic approximation (see, e.g., \cite{li2023online, srikant2024rates, wu2025uncertainty}) through the solution to the Poisson equation and a telescoping technique. We refer readers to Section 2 of \cite{wu2025uncertainty} for a detailed discussion of the mechanism of Markov chain-induced martingales. 


\begin{theorem}\label{thm:Berry-Esseen-Kol-Markov}
Consider a Markov chain $\{s_t\}_{t \geq 0}$ on state space $\mathcal{S}$ with a unique stationary distribution $\mu$ and a spectral gap $1-\lambda > 0$. Let $\bm{f}$ be a function from $\mathcal{S}^2$ to $\mathbb{R}^{d }$ satisfying 
\begin{align*}
&\mathbb{E}_{s' \sim P(\cdot \mid s)}\bm{f}(s,s') = \bm{0}, \quad  \alpha \bm{I} \preceq \mathbb{E}_{s' \sim P(\cdot \mid s)}[\bm{f}(s,s')\bm{f}^\top(s,s')] \preceq \beta\bm{I}, \quad \text{and} \\ 
&\mathbb{E}_{s' \sim P(\cdot \mid s)}\|\bm{f}(s,s')\|_{\infty}^3 \leq (\gamma \log_+d)^{\frac{3}{2}},\quad \forall s \in \mathcal{S},
\end{align*}
for positive constants $\alpha,\beta$ and $\gamma$.
Further define
\begin{align*}
\bm{\Sigma} =\mathbb{E}_{s\sim\mu,s' \sim P(\cdot \mid s)}\Big[\bm{f}(s,s')\bm{f}^\top(s,s')\Big], 
\end{align*} 
and assume that the diagonal elements of $\bm{\Sigma}$ are all equal to 1. \footnote{Since the Kolmogorov distance is invariant to element-wise multiplication, this assumption is not essential to the theorem and only for ease of presentation.} Let $\nu$ denote a probability distribution on $\mathcal{S}$ that is absolutely continuous with respect to the stationary distribution $\mu$; furthermore, assume that there exists $p \in (1,\infty]$, such that the Radon-Nikodym derivative of $\nu$ with respect to $\mu$ satisfies $\left\|\frac{\mathrm{d}\nu}{\mathrm{d}\mu}\right\|_{\mu,p} < \infty$.
Then, when $s_0 \sim \nu$, there exists $\kappa \asymp \log(nd)n^{-\frac{1}{2}}$, such that
\begin{align}
d_{\mathsf{K}}\left(\frac{1}{\sqrt{n}}\sum_{k=1}^n \bm{f}(s_{k-1},s_k),\mathcal{N}(\bm{0},\bm{\Sigma})\right) \lesssim \left(\sqrt{\kappa \log(nd)}+ \|\bm{\Sigma}^{-1}\| \kappa \log \frac{1}{\kappa }\right)\log_+d + \frac{\gamma^{\frac{3}{4}} + \beta^{\frac{3}{4}}}{\sqrt{\alpha}}\log_+^2d \; n^{-\frac{1}{4}} 
\end{align}
\end{theorem}

Observe that Theorem \ref{thm:Berry-Esseen-Kol-Markov} does not depend on Assumption \ref{as:fixed-variation}; instead, it controls the difference between the conditional and unconditional variance through a matrix Hoeffding's inequality on Markov chains, and constructs an auxiliary martingale using Yurinskii's coupling technique \cite{yurinksi.martingale}. A related result in \cite{wu2025uncertainty}, which measures the Gaussian approximation error in convex distance, does not depend on the function $\bm{f}$ being time-invariant and the lower bound on the eigenvalue of the conditional variance.


\begin{section}{Acknowledgments}
    Weichen Wu and Alessandro Rinaldo acknowledge partial support from the NIH Award RF1NS121913. Dung Le acknowledges partial funding from  NSF grant CCF-2505865, supporting the NSF AI Institute for Foundations of Machine Learning (IFML). The authors thank Prof. Yuting Wei from the Wharton School, University of Pennsylvania, for helpful discussions.
\end{section}

\bibliography{main.bib}

@misc{rubtsov2026gaussian,
      title={Gaussian Approximation for Asynchronous {Q}-learning}, 
      author={Artemy Rubtsov and Sergey Samsonov and Vladimir Ulyanov and Alexey Naumov},
      year={2026},
      eprint={2604.07323},
      url={https://arxiv.org/abs/2604.07323}
}

@article{Barbour1990SteinsMF,
  title={{S}tein's method for diffusion approximations},
  author={Andrew D. Barbour},
  journal={Probability Theory and Related Fields},
  year={1990},
  volume={84},
  pages={297--322}
}

@article{meckes2009stein,
  title={On {S}tein's method for multivariate normal approximation},
  author={Meckes, Elizabeth S.},
  journal={arXiv preprint arXiv:0902.0333},
  year={2009}
}

@article{gotze1991rate,
  title={On the rate of convergence in the multivariate {CLT}},
  author={G{\"o}tze, Friedrich},
  journal={The Annals of Probability},
  volume={19},
  number={2},
  pages={724--739},
  year={1991},
  publisher={Institute of Mathematical Statistics},
  doi={10.1214/aop/1176990448}
}

@article{Bhattacharya2010AnEO,
  title={An Exposition of {G}{\"o}tze's Estimation of the Rate of Convergence in the Multivariate Central Limit Theorem},
  author={Rabi Bhattacharya and Susan Holmes},
  journal={arXiv preprint arXiv:1003.4254},
  year={2010}
}

@incollection{Dvoretzky.clt.martingale:72,
author = {Dvoretzky, A.},
booktitle  = {Proceedings of the Sixth Berkeley Symposium on Mathematical Statistics and Probability},
editor      = {Le Cam, L. M. and  Neyman, J. and Scott, E. L.},
pages = {513-535},
title = {Asymptotic normality for sums of dependent random variables},
volume = {6.2},
publisher={University of California Press},
year = {1972},
}

@article{yurinksi.martingale,
author = {Matias D. Cattaneo and Ricardo P. Masini and William G. Underwood},
title = {{Yurinskii’s coupling for martingales}},
volume = {53},
journal = {The Annals of Statistics},
number = {5},
publisher = {Institute of Mathematical Statistics},
pages = {2179 -- 2203},
year = {2025},
doi = {10.1214/25-AOS2538},
URL = {https://doi.org/10.1214/25-AOS2538}
}

@article{Haausler:88,
author = {Erich Haeusler},
title = {{On the Rate of Convergence in the Central Limit Theorem for Martingales with Discrete and Continuous Time}},
volume = {16},
journal = {The Annals of Probability},
number = {1},
publisher = {Institute of Mathematical Statistics},
pages = {275 -- 299},
year = {1988},
doi = {10.1214/aop/1176991901},
URL = {https://doi.org/10.1214/aop/1176991901}
}

@phdthesis{grams1972rates,
  title  = {Rates of Convergence in the Central Limit Theorem for Dependent Variables},
  author = {Grams, William F.},
  year   = {1972},
  school = {Florida State University},
  address = {Tallahassee, FL}
}

@article{kato1979rates,
  title={Rates of convergence in Central Limit Theorem for martingale differences},
  author={Kato, Yutaka},
  journal={Bulletin of Mathematical Statistics},
  volume={18},
  number={1-2},
  pages={1--8},
  year={1979},
  publisher={Research Association of Statistical Sciences}
}

@article{wu2020berry,
  title={A {B}erry-{E}sseen bound of order $1/\sqrt{n}$ for martingales},
  author={Wu, Songqi and Ma, Xiaohui and Sang, Hailin and Fan, Xiequan},
  journal={Comptes Rendus. Math{\'e}matique},
  volume={358},
  number={6},
  pages={701--712},
  year={2020},
  publisher={Acad{\'e}mie des sciences, Paris}
}

@article{renz1996note,
  title={A note on exact convergence rates in some martingale Central Limit Theorems},
  author={Renz, Joachim},
  journal={The Annals of Probability},
  volume={24},
  number={3},
  pages={1616--1637},
  year={1996},
  publisher={Institute of Mathematical Statistics}
}

@article{mourrat2013rate,
  title     = {On the rate of convergence in the martingale Central Limit Theorem},
  author    = {Mourrat, Jean-Christophe},
  journal   = {Bernoulli},
  volume    = {19},
  number    = {2},
  pages     = {633--645},
  year      = {2013},
  publisher = {Bernoulli Society for Mathematical Statistics and Probability},
  doi       = {10.3150/12-BEJ417},
  url       = {https://projecteuclid.org}
}

@article{Fan2019Exact,
  author = {Fan, Xiequan},
  title = {Exact rates of convergence in some martingale Central Limit Theorems},
  journal = {Journal of Mathematical Analysis and Applications},
  volume = {469},
  number = {2},
  pages = {1028-1044},
  year = {2019},
  doi = {10.1016/j.jmaa.2018.09.049},
  url = {https://doi.org/10.1016/j.jmaa.2018.09.049}
}

@article{rinott2000bounds,
  title={Some bounds on the rate of convergence in the {CLT} for martingales. {II}},
  author={Rinott, Yosef and Rotar, Vladimir I.},
  journal={Theory of Probability \& Its Applications},
  volume={44},
  number={3},
  pages={523--536},
  year={2000},
  publisher={SIAM}
}

@article{el2007exact,
  title={Exact convergence rates in the Central Limit Theorem for a class of martingales},
  author={El Machkouri, Mohamed and Ouchti, Lounes},
  journal={Bernoulli},
  volume={13},
  number={4},
  pages={981--999},
  year={2007},
  publisher={Bernoulli Society for Mathematical Statistics and Probability},
  doi={10.3150/07-BEJ6116},
  url={https://projecteuclid.org/journals/bernoulli/volume-13/issue-4/Exact-convergence-rates-in-the-central-limit-theorem-for-a/10.3150/07-BEJ6116.full}
}

@book{hall2014martingale,
  title     = {Martingale Limit Theory and Its Application},
  author    = {Hall, Peter and Heyde, Christopher C.},
  year      = {1980},
  isbn      = {978-0-12-319350-6}
}

@book{nourdin.peccati.2012,
  publisher = {Cambridge University Press},
  title = {Normal Approximations with Malliavin Calculus.
From Stein's Method to Universality},
  year = {2012},
author = {Ivan Nourdin and Giovanni Peccati}
}

@misc{belloni2018highdimensionalcentrallimit,
      title={A high dimensional Central Limit Theorem for martingales, with applications to context tree models}, 
      author={Alexandre Belloni and Roberto I. Oliveira},
      year={2018},
      eprint={1809.02741},
      archivePrefix={arXiv},
      primaryClass={math.ST},
      url={https://arxiv.org/abs/1809.02741}, 
}

@misc{lopes2022CLT,
      title={Central Limit Theorem and Bootstrap Approximation in High Dimensions: Near $1/\sqrt{n}$ Rates via Implicit Smoothing}, 
      author={Miles E. Lopes},
      year={2022},
      eprint={2009.06004},
      archivePrefix={arXiv},
      primaryClass={math.ST},
      url={https://arxiv.org/abs/2009.06004}, 
}

@article{neeman2024concentration,
    author = {Neeman, Joe and Shi, Bobby and Ward, Rachel},
    title = {Concentration inequalities for sums of {M}arkov-dependent random matrices},
    journal = {Information and Inference: A Journal of the IMA},
    volume = {13},
    number = {4},
    pages = {iaae032},
    year = {2024},
    month = {12},
    issn = {2049-8772},
    doi = {10.1093/imaiai/iaae032},
    url = {https://doi.org/10.1093/imaiai/iaae032},
    eprint = {https://academic.oup.com/imaiai/article-pdf/13/4/iaae032/60957939/iaae032.pdf},
}

@misc{gallouet2018regularity,
      title={Regularity of solutions of the Stein equation and rates in the multivariate Central Limit Theorem}, 
      author={Thomas Gallouët and Guillaume Mijoule and Yvik Swan},
      year={2018},
      eprint={1805.01720},
      archivePrefix={arXiv},
      primaryClass={math.PR},
      url={https://arxiv.org/abs/1805.01720}, 
}

@misc{JMLR2019CLT,
  title={Normal Approximation for {S}tochastic {G}radient {D}escent via Non-Asymptotic Rates of Martingale {CLT}},
  author={Anastasiou, Andreas and Balasubramanian, Krishnakumar and Erdogdu, Murat A.},
  journal = {Proceedings of the Thirty-Second Conference on Learning Theory},
  pages = {115--137},
  year = {2019},
  editor = 	 {Beygelzimer, Alina and Hsu, Daniel},
  volume = 	 {99},
  series = 	 {Proceedings of Machine Learning Research},
  publisher =    {PMLR},
  url = 	 {https://proceedings.mlr.press/v99/anastasiou19a.html},
}

@article{rollin2018,
   title={On quantitative bounds in the mean martingale Central Limit Theorem},
   volume={138},
   ISSN={0167-7152},
   url={http://dx.doi.org/10.1016/j.spl.2018.03.004},
   DOI={10.1016/j.spl.2018.03.004},
   journal={Statistics \& Probability Letters},
   publisher={Elsevier BV},
   author={Röllin, Adrian},
   year={2018},
   month=jul, pages={171–176} }

@misc{arun2020hdCLT,
      title={High-dimensional {CLT} for Sums of Non-degenerate Random Vectors: $n^{-1/2}$-rate}, 
      author={Arun Kumar Kuchibhotla and Alessandro Rinaldo},
      year={2020},
      eprint={2009.13673},
      archivePrefix={arXiv},
      primaryClass={math.ST},
      url={https://arxiv.org/abs/2009.13673}, 
}

@misc{nourdin2021multivariate,
      title={Multivariate normal approximation on the {W}iener space: new bounds in the convex distance}, 
      author={Ivan Nourdin and Giovanni Peccati and Xiaochuan Yang},
      year={2021},
      eprint={2001.02188},
      archivePrefix={arXiv},
      primaryClass={math.PR},
      url={https://arxiv.org/abs/2001.02188}, 
}

@article{Kojevnikov2022,
   title={A {B}erry–{E}sseen bound for vector-valued martingales},
   volume={186},
   ISSN={0167-7152},
   url={http://dx.doi.org/10.1016/j.spl.2022.109448},
   DOI={10.1016/j.spl.2022.109448},
   journal={Statistics \& Probability Letters},
   publisher={Elsevier BV},
   author={Kojevnikov, Denis and Song, Kyungchul},
   year={2022},
   pages={109448} 
}

@misc{srikant2024rates,
      title={Rates of Convergence in the Central Limit Theorem for {Markov} Chains, with an Application to {TD} learning}, 
      author={R. Srikant},
      year={2026},
      eprint={2401.15719},
      archivePrefix={arXiv},
      primaryClass={math.PR},
      url={https://arxiv.org/abs/2401.15719}, 
}

@article{li2023online,
  title={Online Statistical Inference for Nonlinear Stochastic Approximation with Markovian Data},
  author={Li, Xiang and Liang, Jiadong and Zhang, Zhihua},
  journal={arXiv preprint arXiv:2302.07690},
  year={2023}
}

@article{chernozhukov2017central,
  author  = {Chernozhukov, Victor and Chetverikov, Denis and Kato, Kengo},
  title   = {Central Limit Theorems and Bootstrap in High Dimensions},
  journal = {The Annals of Probability},
  volume  = {45},
  number  = {4},
  pages   = {2309--2352},
  year    = {2017},
  doi     = {10.1214/16-AOP1117}
}

@article{fang2021highdimensional,
    title = {High-dimensional Central Limit Theorems by {S}tein's method},
    author = {Fang, Xiao and Koike, Yuta},
    journal = {Annals of Applied Probability},
    volume = {31},
    issue = {4},
    pages = {1660-1686},
    year = {2021}
}

@misc{wu2025uncertainty,
      title={Uncertainty quantification for {M}arkov chains with application to {T}emporal {D}ifference learning}, 
      author={Weichen Wu and Yuting Wei and Alessandro Rinaldo},
      year={2025},
      eprint={2502.13822},
      archivePrefix={arXiv},
      primaryClass={stat.ML},
      url={https://arxiv.org/abs/2502.13822}, 
}

@article{bothausen1982convergence,
author = {E. Bolthausen},
title = {{Exact Convergence Rates in Some Martingale Central Limit Theorems}},
volume = {10},
journal = {The Annals of Probability},
number = {3},
publisher = {Institute of Mathematical Statistics},
pages = {672 -- 688},
year = {1982},
doi = {10.1214/aop/1176993776},
URL = {https://doi.org/10.1214/aop/1176993776}
}
\bibliographystyle{plainnat}

\appendix 

\section{Preliminary facts}
\begin{theorem}\cite[Lemma 2.3]{fang2021highdimensional}\label{thm:Fang}
For every $\mathcal{A} \in \mathscr{H}_d$, $\sigma>0$ and $\bm{x} \in \mathbb{R}^d$, define
\begin{align}\label{eq:defn-phi}
\varphi_{\mathcal{A},\sigma}(\bm{x}) = \mathbb{P}_{\eta \sim \mathcal{N}(\bm{0},\bm{I}_d)}(\bm{x} + \sigma\bm{\eta} \in \mathcal{A}).
\end{align}
Then it can be guaranteed that for every $s \in \mathbb{N}_+$, there exists a constant $C_s$ depending only on $s$, such that
\begin{align*}
\|\nabla^s \varphi_{\mathcal{A},\sigma}(\bm{x})\|_1 \leq C_s\sigma^{-s}\log_+^{\frac{s}{2}} d.
\end{align*}
\end{theorem}


\begin{theorem}\cite[Lemma 2.3]{Kojevnikov2022}\label{thm:Gaussian-inf-norm}
Let $\bm{X}$ denote a $d$-dimensional Gaussian random variable with mean $\bm{0}$ and variance $\bm{\Sigma}$. Then for every $s \geq 2$, there exists a constant $C_s$ depending only on $s$, such that
\begin{align*}
\mathbb{E}\|\bm{X}\|_{\infty}^s \leq C_s \log_+^{\frac{s}{2}}d \cdot (d_{\max}(\bm{\Sigma}))^{\frac{s}{2}};
\end{align*}
furthermore, for any $\delta \in (0,1)$, 
\begin{align*}
\|\bm{X}\|_{\infty} \lesssim \sqrt{d_{\max}(\bm{\Sigma}) \log \frac{d}{\delta}}.
\end{align*}
with probability at least $1-\delta$.
\end{theorem}

\begin{theorem}[Gaussian comparison in Kolmogorov distance]\cite[Theorem 2.3]{lopes2022CLT}\label{thm:Gaussian-comparison-Kol}
Let $\bm{Z} \sim \mathcal{N}(\bm{0},\bm{\Sigma})$ and $\bm{Z}' \sim \mathcal{N}(\bm{0},\bm{\Sigma}')$ where $\bm{\Sigma}$ and $\bm{\Sigma}'$ are $d \times d$ covariance matrices. Denote $\bm{D} = \mathsf{diag}\{\Sigma_{11},\Sigma_{22},...,\Sigma_{dd}\}$, and define
\begin{align*}
\Delta = \|\bm{D}^{-\frac{1}{2}}(\bm{\Sigma}' - \bm{\Sigma})\bm{D}^{-\frac{1}{2}}\|_{e,\infty}.
\end{align*}
Then the Kolmogorov distance between $\bm{Z}$ and $\bm{Z}'$ is bounded by
\begin{align*}
d_{\mathsf{K}}(\bm{Z},\bm{Z}') \lesssim \|\bm{D}^{\frac{1}{2}}\bm{\Sigma}^{-1}\bm{D}^{\frac{1}{2}}\| \Delta \log \frac{1}{\Delta}\log_+d.
\end{align*}
\end{theorem}

\begin{theorem}[Gaussian anti-concentration]\cite[Lemma A.1]{chernozhukov2017central}\label{thm:anti-concentration}
Let $\bm{Z} \sim \mathcal{N}(\bm{0},\bm{\Sigma})$ be a $d$-dimensional Gaussian random variable and $\mathcal{A} \in \mathscr{H}_d$. Define
\begin{align*}
\mathcal{A}^{\varepsilon} = \{\bm{x} \in \mathbb{R}^d: \inf_{y \in \mathcal{A}}\|\bm{x} - \bm{y}\|_{\infty} \leq \varepsilon\}, \quad \text{and} \quad \mathcal{A}^{-\varepsilon} = \{\bm{x} \in \mathcal{A}: \forall \bm{y} \text{ with } \|\bm{x} - \bm{y}\|_{\infty} \leq \varepsilon, \bm{y} \in \mathcal{A}\}.
\end{align*}
Then it can be guaranteed that
\begin{align*}
\mathbb{P}(\bm{Z} \in \mathcal{A}^\varepsilon) - \mathbb{P}(\bm{Z} \in \mathcal{A}) \lesssim \sqrt{\frac{\log_+d}{d_{\min}(\bm{\Sigma})}}\varepsilon,\quad \text{and} \quad \mathbb{P}(\bm{Z} \in \mathcal{A}) - \mathbb{P}(\bm{Z} \in \mathcal{A}^{-\varepsilon}) \lesssim \sqrt{\frac{\log_+d}{d_{\min}(\bm{\Sigma})}}\varepsilon.
\end{align*}
\end{theorem}

\section{Proofs}

\subsection{Proof of Theorems \ref{thm:Berry-Esseen-Kol-naive} and \ref{thm:Berry-Esseen-Kol}}\label{app:proof-thm-Berry-Esseen-Kol}



Throughout the proof, we let
\begin{align*}
\bm{S}_k = \sum_{j=1}^k \bm{X}_j, \quad \bm{Y}_j = \bm{V}_j^{\frac{1}{2}}\bm{Z}_j, \quad \bm{T}_k = \sum_{j=k}^n \bm{Y}_j,
\end{align*}
where $\{\bm{Z}_j\}_{j \in [n]}$ denote a sequence of $i.i.d.$ standard Gaussian random variables in $\mathbb{R}^d$ that are independent of $\{\bm{X}_k\}_{k \in [n]}$. 
Furthermore, let 
\begin{align*}
\bm{P}_k = \sum_{j=k}^n \bm{V}_j
\end{align*}
and recall  that, by assumption, $\bm{P}_1 = n \bm{\Sigma}_n$.

Following the same arguments as in the beginning of Theorem 1 of \cite{fang2021highdimensional} (see, in particular, Eqs. (2.4) and (2.6) therein), we obtain that for any $t>0$,
\begin{align}\label{eq:Berry-Esseen-Kol-decompose-1}
&d_{\mathsf{K}}\left(\frac{1}{\sqrt{n}}\sum_{k=1}^n \bm{X}_k,\mathcal{N}(\bm{0},\bm{\Sigma}_n)\right)\\
&=d_{\mathsf{K}}(\bm{S}_n,\bm{T}_1) \nonumber \\ 
&\lesssim d_{\mathsf{K}}(e^{-t}\bm{S}_n + \sqrt{1-e^{-2t}}\bm{T}_1,\bm{T}_1) + \frac{\sqrt{1-e^{-2t}}}{e^{-t}}\sqrt{\frac{d_{\max}(\bm{\Sigma}_n)}{d_{\min}(\bm{\Sigma}_n)}} \log_+d \nonumber \\ 
&= \sup_{\mathcal{A} \in \mathscr{H}_d}\left|\mathbb{P}(e^{-t}\bm{S}_n + \sqrt{1-e^{-2t}}\bm{P}_1^{\frac{1}{2}}\bm{Z} \in \mathcal{A}) - \mathbb{P}(\bm{P}_1^{\frac{1}{2}}\bm{Z} \in \mathcal{A})\right|+\frac{\sqrt{1-e^{-2t}}}{e^{-t}}\sqrt{\frac{d_{\max}(\bm{\Sigma}_n)}{d_{\min}(\bm{\Sigma}_n)}} \log_+d
\end{align}

We will specify a particular value of $t$ at the end of the proof. Notice that on the third line, we invoked the fact that $\bm{T}_1 \sim \mathcal{N}(\bm{0},\bm{P}_1)$. Meanwhile, the basic property of Gaussian distributions implies
\begin{align*}
\bm{T}_1 \overset{d}{=} e^{-t}\bm{T}_1 + \sqrt{1-e^{-2t}}\bm{P}_1^{\frac{1}{2}}\bm{Z},
\end{align*}
where $\bm{Z}$ is a $d$-dimensional random vector with a standard Gaussian distribution, independent of  $\{\bm{Z}_j\}_{j \in [n]}$ and of $\{\bm{X}_k\}_{k \in [n]}$.
Therefore, by defining for each $\mathcal{A} \in \mathscr{H}_d$, the (smooth) function
\begin{align*}
\bm{x} \in \mathbb{R}^d \mapsto h_{\mathcal{A},t}(\bm{x}) = \mathbb{P}(e^{-t}\bm{x}+ \sqrt{1-e^{-2t}}\bm{P}_1^{\frac{1}{2}}\bm{Z} \in \mathcal{A}),
\end{align*}
we can further simplify \eqref{eq:Berry-Esseen-Kol-decompose-1} as
\begin{align}\label{eq:Berry-Esseen-Kol-decompose-2}
d_{\mathsf{K}}(\bm{S}_n,\bm{T}_1) &\lesssim \sup_{\mathcal{A} \in \mathscr{H}_d}|\mathbb{E}[h_{\mathcal{A},t}(\bm{S}_n)-h_{\mathcal{A},t}(\bm{T}_1)]| +\frac{\sqrt{1-e^{-2t}}}{e^{-t}}\sqrt{\frac{d_{\max}(\bm{\Sigma}_n)}{d_{\min}(\bm{\Sigma}_n)}} \log_+d.
\end{align}

For the first term on the right-hand-side, a classic Lindeberg swapping argument yields
\begin{align}\label{eq:Lindeberg-swapping}
&\sup_{\mathcal{A} \in \mathscr{H}_d}|\mathbb{E}[h_{\mathcal{A},t}(\bm{S}_n)-h_{\mathcal{A},t}(\bm{T}_1)]|\nonumber \\ 
&\leq \sup_{\mathcal{A} \in \mathscr{H}_d} \left|\sum_{k=1}^n \mathbb{E}\left[h_{\mathcal{A},t}(\bm{S}_{k-1} + \bm{X}_k + \bm{T}_{k+1})-h_{\mathcal{A},t}(\bm{S}_{k-1} + \bm{Y}_k + \bm{T}_{k+1})\right]\right| \nonumber \\ 
&\leq  \sup_{\mathcal{A} \in \mathscr{H}_d} \sum_{k=1}^n \mathbb{E} \left[ \left|\mathbb{E}[h_{\mathcal{A},t}(\bm{S}_{k-1} + \bm{X}_k + \bm{T}_{k+1})-h_{\mathcal{A},t}(\bm{S}_{k-1} + \bm{Y}_k + \bm{T}_{k+1})\mid \mathscr{F}_{k-1}] \right| \right],
\end{align}
where we set $\bm{S}_{0} = \bm{T}_{n+1} = \bm{0}$ and in the third line we have used the triangle inequality, the law of iterated expectations and the convexity of the $\sup$ function.

Notice that $\bm{S}_{k-1}$ is measurable with respect to $\mathscr{F}_{k-1}$, and for every hyper-rectangle $\mathcal{A}$, the set
\begin{align*}
\{\bm{x}-\bm{S}_{k-1}: \bm{x} \in \mathcal{A} \}
\end{align*}
is also a hyper-rectangle. Meanwhile, due to Assumption \ref{as:fixed-variation}, the matrices $\bm{P}_k$ and $\bm{P}_{k+1}$ are both measurable with respect to $\mathscr{F}_{k-1}$. Consequently, the inner conditional expectation of each summand in \eqref{eq:Lindeberg-swapping} can  be written as
 \begin{align}\label{eq:Lindeberg-summand}
& \mathbb{E}[h_{\mathcal{A},t}(\bm{S}_{k-1} + \bm{X}_k + \bm{T}_{k+1})-h_{\mathcal{A},t}(\bm{S}_{k-1} + \bm{Y}_k + \bm{T}_{k+1})\mid \mathscr{F}_{k-1}] \nonumber \\ 
&= \mathbb{P}(e^{-t}(\bm{X}_k + \bm{T}_{k+1}) + \sqrt{1-e^{-2t}}\bm{P}_1^{\frac{1}{2}}\bm{Z} \in \mathcal{A}\mid \mathscr{F}_{k-1})\nonumber \\
&\qquad \qquad - \mathbb{P}(e^{-t}(\bm{Y}_k + \bm{T}_{k+1}) + \sqrt{1-e^{-2t}}\bm{P}_1^{\frac{1}{2}}\bm{Z} \in \mathcal{A}\mid \mathscr{F}_{k-1}).
 \end{align}

 For the second term, observe that
 \begin{align*}
\bm{Y}_k + \bm{T}_{k+1}\mid \mathscr{F}_{k-1}\sim \mathcal{N}(\bm{0},\bm{P}_k),
 \end{align*}
 so
 \begin{align*}
e^{-t}(\bm{Y}_k + \bm{T}_{k+1}) + \sqrt{1-e^{-2t}}\bm{P}_1^{\frac{1}{2}}\bm{Z} \mid \mathscr{F}_{k-1} \overset{d}{=} \bm{P}_k^{\frac{1}{2}}\bm{Z} + \sqrt{1-e^{-2t}}(\bm{P}_1-\bm{P}_k)^{\frac{1}{2}}\bm{Z}'
 \end{align*}
where $\bm{Z}'$ is another standard Gaussian random variable, independent of $\bm{Z}$, $\{\bm{Z}_j\}_{j \in [n]}$ and $\{\bm{X}_k\}_{k \in [n]}$. Similarly, we can write
\begin{align*}
&e^{-t}(\bm{X}_k + \bm{T}_{k+1}) + \sqrt{1-e^{-2t}}\bm{P}_1^{\frac{1}{2}}\bm{Z} \mid \mathscr{F}_{k-1} \\ 
&\overset{d}{=}e^{-t}(\bm{X}_k + \bm{T}_{k+1}) + \sqrt{1-e^{-2t}} \bm{P}_k^{\frac{1}{2}}\bm{Z} + \sqrt{1-e^{-2t}}(\bm{P}_1-\bm{P}_k)^{\frac{1}{2}}\bm{Z}'.
\end{align*}
Hence by denoting
\begin{align*}
\bm{x} \in \mathbb{R}^D \mapsto \varphi_k(\bm{x}) := \mathbb{P}(\bm{x} + \sqrt{1-e^{-2t}}(\bm{P}_1-\bm{P}_k)^{\frac{1}{2}}\bm{Z}' \in \mathcal{A} | \mathcal{F}_{k-1} ),
\end{align*}
a random function depending on $\bm{P}_k$ and $\mathcal{A}$, the Lindeberg summand \eqref{eq:Lindeberg-summand} can be written as 
\begin{align}\label{eq:Lindeberg-summand-2}
&\mathbb{P}(e^{-t}(\bm{X}_k + \bm{T}_{k+1}) + \sqrt{1-e^{-2t}}\bm{P}_1^{\frac{1}{2}}\bm{Z} \in \mathcal{A}\mid \mathscr{F}_{k-1})-\mathbb{P}(e^{-t}(\bm{Y}_k + \bm{T}_{k+1}) + \sqrt{1-e^{-2t}}\bm{P}_1^{\frac{1}{2}}\bm{Z} \in \mathcal{A}\mid \mathscr{F}_{k-1})\nonumber \\ 
&= \mathbb{E}_{k-1}[\varphi_k(e^{-t}(\bm{X}_k + \bm{T}_{k+1}) + \sqrt{1-e^{-2t}}\bm{P}_k^{\frac{1}{2}}\bm{Z})] -\mathbb{E}_{k-1}[\varphi_k(\bm{P}_k^{\frac{1}{2}}\bm{Z})].
\end{align}
Next, we apply Stein's method, conditionally on $\mathcal{F}_{k-1}$. In detail, for any $\bm{x} \in \mathbb{R}^d$, consider the random function
\begin{align}\label{eq:defn-psi}
\bm{x} \in \mathbb{R}^D \mapsto  \psi_k(\bm{x}):= -\int_t^{\infty} \mathbb{E}_{k-1}\left[\varphi_k(e^{-s}\bm{x} + \sqrt{1-e^{-2s}}\bm{P}_k^{\frac{1}{2}}\bm{Z})\right] - \mathbb{E}_{k-1}\left[\varphi_k (\bm{P}_k^{\frac{1}{2}}\bm{Z})\right]\mathrm{d}s.
\end{align}


Then, $\psi_k$ is the solution to the conditional Stein's equation
\begin{align}\label{eq:Stein}
\langle \bm{P}_k,\bm\nabla^2 \psi_k(\bm{x})\rangle - \langle \bm{x},\bm\nabla \psi_k(\bm{x})\rangle = \mathbb{E}_{k-1}[\varphi_k(e^{-t}\bm{x} + \sqrt{1-e^{-2t}}\bm{P}_k^{\frac{1}{2}}\bm{Z})] - \mathbb{E}_{k-1}[\varphi_k(\bm{P}_k^{\frac{1}{2}}\bm{Z})], \quad \forall \bm{x} \in \mathbb{R}^d.
\end{align}
See, e.g., \cite[Eq. 1.16 in][]{Bhattacharya2010AnEO} and \cite{gotze1991rate}.
Thus, the Lindeberg summand \eqref{eq:Lindeberg-summand-2} can be further decomposed as
\begin{align}\label{eq:Lindeberg-summand-3}
&\mathbb{E}_{k-1}[\varphi_k(e^{-t}(\bm{X}_k + \bm{T}_{k+1}) + \sqrt{1-e^{-2t}}\bm{P}_k^{\frac{1}{2}}\bm{Z})] -\mathbb{E}_{k-1}[\varphi_k(\bm{P}_k^{\frac{1}{2}}\bm{Z})]\nonumber \\ 
&= \mathbb{E}_{k-1}[\langle\bm{P}_k,\bm\nabla^2 \psi_k(\bm{X}_k + \bm{T}_{k+1})\rangle - \langle \bm{X}_k + \bm{T}_{k+1},\bm{\nabla}\psi_k(\bm{X}_k + \bm{T}_{k+1})\rangle]\nonumber \\ 
&\overset{(i)}{=} \mathbb{E}_{k-1}[\langle\bm{P}_k,\bm\nabla^2 \psi_k(\bm{X}_k + \bm{T}_{k+1})\rangle] - \mathbb{E}_{k-1}[\langle \bm{X}_k,\bm{\nabla}\psi_k(\bm{X}_k + \bm{T}_{k+1})\rangle]-\mathbb{E}_{k-1}[\langle\bm{P}_{k+1},\bm\nabla^2 \psi_k(\bm{X}_k + \bm{T}_{k+1})\rangle]\nonumber \\ 
&= \mathbb{E}_{k-1}[\langle\bm{V}_k,\bm\nabla^2 \psi_k(\bm{X}_k + \bm{T}_{k+1})\rangle] - \mathbb{E}_{k-1}[\langle \bm{X}_k,\bm{\nabla}\psi_k(\bm{X}_k + \bm{T}_{k+1})\rangle] \nonumber \\ 
&\overset{(ii)}{=} \mathbb{E}_{k-1}[\langle\bm{V}_k,\bm\nabla^2 \psi_k( \bm{T}_{k+1})\rangle] + \mathbb{E}_{k-1}\left[\int_{\theta=0}^1\langle\bm\nabla^3 \psi_k(\theta\bm{X}_k + \bm{T}_{k+1}),\bm{V}_k \otimes \bm{X}_k \rangle \mathrm{d}\theta \right] \nonumber \\ 
&- \mathbb{E}_{k-1}[\langle \bm{X}_k, \bm\nabla\psi_k(\bm{T}_{k+1})\rangle] - \mathbb{E}_{k-1}[\langle \bm{X}_k^{\otimes 2}, \bm\nabla^2\psi_k(\bm{T}_{k+1})\rangle]-\frac{1}{2}\mathbb{E}_{k-1}\left[\int_{\theta=0}^1\langle \bm{X}_k^{\otimes 3}, \bm\nabla^3\psi_k(\theta\bm{X}_k+\bm{T}_{k+1})\rangle \mathrm{d}\theta\right].
\end{align}
Above, we invoked Stein's equation for non-isotropic centered Gaussian \citep[see, e.g.][]{nourdin.peccati.2012} in (i), and Taylor's expansion in (ii). 
Furthermore, since $\mathbb{E}_{k-1}[\bm{X}_k] = \bm{0}$ and $\mathbb{E}_{k-1}[\bm{X}_k^{\otimes 2}]=\bm{V}_k$, Equation \eqref{eq:Lindeberg-summand-3} can be further simplified as
\begin{align*}
&\mathbb{E}_{k-1}[\varphi_k(e^{-t}(\bm{X}_k + \bm{T}_{k+1}) + \sqrt{1-e^{-2t}}\bm{P}_k^{\frac{1}{2}}\bm{Z})] -\mathbb{E}_{k-1}[\varphi_k(\bm{P}_k^{\frac{1}{2}}\bm{Z})]\nonumber \\ 
&=  \mathbb{E}_{k-1}\left[\int_{\theta=0}^1\langle\bm\nabla^3 \psi_k(\theta\bm{X}_k + \bm{T}_{k+1}),\bm{V}_k \otimes \bm{X}_k \rangle \mathrm{d}\theta \right]-\frac{1}{2}\mathbb{E}_{k-1}\left[\int_{\theta=0}^1\langle \bm{X}_k^{\otimes 3}, \bm\nabla^3\psi_k(\theta\bm{X}_k+\bm{T}_{k+1})\rangle \mathrm{d}\theta\right],
\end{align*}
while the triangle inequality and Hölder's inequality yield
\begin{align}\label{eq:Lindeberg-summand-4}
&\left|\mathbb{E}_{k-1}[\varphi_k(e^{-t}(\bm{X}_k + \bm{T}_{k+1}) + \sqrt{1-e^{-2t}}\bm{P}_k^{\frac{1}{2}}\bm{Z})] -\mathbb{E}_{k-1}[\varphi_k(\bm{P}_k^{\frac{1}{2}}\bm{Z})]\right|\nonumber \\
&\leq 2\sup_{\bm{x} \in \mathbb{R}^d}\|\bm\nabla^3\psi_k(\bm{x})\|_1 \mathbb{E}_{k-1}\|\bm{X}_k\|_{\infty}^3.
\end{align}
Next, we proceed to bound the term
\[
\sup_{\bm{x} \in \mathbb{R}^d}\|\bm\nabla^3\psi_k(\bm{x})\|_1.
\]
First we note that, using the definition of $\psi_k$ in \eqref{eq:defn-psi} and by integrability, for any $\bm{x} \in \mathbb{R}^d$,
\begin{align}\label{eq:nabla3-psi-k}
\bm\nabla^3 \psi_k(\bm{x}) &= -\int_t^\infty \bm{\nabla}^3 \varphi_k(e^{-s}\bm{x} + \sqrt{1-e^{-2s}}\bm{P}_k^{\frac{1}{2}}\bm{Z}) \mathrm{d}s \nonumber \\ 
&= -\int_t^\infty \bm{\nabla}^3 \mathbb{P}(e^{-s}\bm{x} + \sqrt{1-e^{-2s}}\bm{P}_k^{\frac{1}{2}}\bm{Z} +  \sqrt{1-e^{-2t}}(\bm{P}_1-\bm{P}_k)^{\frac{1}{2}}\bm{Z}' \in \mathcal{A} \mid \mathcal{F}_{k-1})\mathrm{d}s.
\end{align}


Since $\bm{Z}$ and $\bm{Z}'$ are independent, it holds that, for every $s \geq t$, 
\begin{align*}
\sqrt{1-e^{-2s}}\bm{P}_k^{\frac{1}{2}}\bm{Z} +  \sqrt{1-e^{-2t}}(\bm{P}_1-\bm{P}_k)^{\frac{1}{2}}\bm{Z}' | \mathcal{F}_{k-1} \; \sim \mathcal{N}(\bm{0},(1-e^{2t})\bm{P}_1+(e^{-2t}-e^{-2s})\bm{P}_k) | \mathcal{F}_{k-1}
\end{align*}
For simplicity, set 
\begin{align*}
\varepsilon_k := \lambda_{\min}(\bm{P}_k)
\end{align*}
and for a fixed $s$,
\begin{align*}
\bm{\Sigma}_{t,k} := (1-e^{2t})\bm{P}_1+(e^{-2t}-e^{-2s})\bm{P}_k,
\end{align*}
both random quantities measurable with respect to $\mathcal{F}_{k-1}$.
Then, almost surely, 
\begin{align*}
\varepsilon_{t,k} := (1-e^{2t})\varepsilon_1 + (e^{-2t}-e^{-2s})\varepsilon_k \leq \lambda_{\min}(\bm{\Sigma}_{t,k}).
\end{align*}
Next, by independence and using the properties of the Gaussian distribution, 
\[
\mathbb{P}(e^{-s}\bm{x} + \sqrt{1-e^{-2s}}\bm{P}_k^{\frac{1}{2}}\bm{Z} +  \sqrt{1-e^{-2t}}(\bm{P}_1-\bm{P}_k)^{\frac{1}{2}}\bm{Z}' \in \mathcal{A} | \mathcal{F}_{k-1}) 
\]
can be equivalently expressed as
\begin{align*}
\mathbb{P}(e^{-s}\bm{x} + \Sigma_{k,t}^{1/2} \bm{Z}'' \in \mathcal{A} | \mathcal{F}_{k-1}) & = \mathbb{P}(e^{-s}\bm{x} + (\bm{\Sigma}_{t,k} - \varepsilon_{t,k} \bm{I}_d)^{1/2}  \bm{Z}'' + \sqrt{\varepsilon_{t,k}} \bm{I}_d \bm{Z}'''\in \mathcal{A} | \mathcal{F}_{k-1})\\ 
& \leq \sup_{ \bm{y} \in \mathbb{R}^d} \mathbb{P}(\bm{y} + \sqrt{\varepsilon_{t,k}} \bm{I}_d \bm{Z}''' \in \mathcal{A} | \mathcal{F}_{k-1} ),
\end{align*}
where $\bm{Z}''$ and $\bm{Z}'''$ are independent $d$-dimensional standard Gaussian vectors, also independent of $\mathcal{F}_{k-1}$. Theorem \ref{thm:Fang} guarantees that for every $\bm{x} \in \mathbb{R}^d$ and $\mathcal{A} \in \mathscr{H}_d$,
\begin{align*} 
&\left\|\bm{\nabla}^3 \mathbb{P}(e^{-s}\bm{x} + \sqrt{1-e^{-2s}}\bm{P}_k^{\frac{1}{2}}\bm{Z} +  \sqrt{1-e^{-2t}}(\bm{P}_1-\bm{P}_k)^{\frac{1}{2}}\bm{Z}' \in \mathcal{A} | \mathcal{F}_{k-1})\right\|_1 \\ 
& \leq \left\| \sup_{ \bm{y} \in \mathbb{R}^d} \mathbb{P}(\bm{y} + \sqrt{\varepsilon_{t,k}} \bm{I}_d \bm{Z}''' \in \mathcal{A} | \mathcal{F}_{k-1} ) \right\|_1 \\
& \lesssim  \frac{e^{-3s}}{ \varepsilon_{t,k}^{\frac{3}{2}}} \cdot \log_+^{\frac{3}{2}}d\\
&= \frac{e^{-3s}}{((1-e^{2t})\varepsilon_1 + (e^{-2t}-e^{-2s})\varepsilon_k)^{\frac{3}{2}}} \cdot \log_+^{\frac{3}{2}}d,
\end{align*}
almost surely.
Plugging this bound into \eqref{eq:nabla3-psi-k}, we obtain, for every $\bm{x} \in \mathbb{R}^d$ and $\mathcal{A} \in \mathscr{H}_d$,
\begin{align}\label{eq:nabla3-psi-k-bound-1}
\|\bm\nabla^3 \psi_k(\bm{x})\|_1 \lesssim \log_+^{\frac{3}{2}}d  \cdot \int_t^\infty\frac{e^{-3s}}{((1-e^{2t})\varepsilon_1 + (e^{-2t}-e^{-2s})\varepsilon_k)^{\frac{3}{2}}} \mathrm{d}s,
\end{align}
almost surely.
The integral on the right hand side of the previous expression can be conveniently bounded almost surely as
\begin{align}\label{eq:prop.b1}
&\int_t^\infty\frac{e^{-3s}}{((1-e^{2t})\varepsilon_1 + (e^{-2t}-e^{-2s})\varepsilon_k)^{\frac{3}{2}}} \mathrm{d}s \leq \frac{e^{-3t}}{\sqrt{(1-e^{-2t})\varepsilon_1}} \cdot \frac{1}{(1-e^{-2t})\varepsilon_1 + e^{-2t}\varepsilon_k}.
\end{align}
For a proof of the previous inequality, see Section \ref{app:proof-Stein-integral}. 
Plugging this bound into \eqref{eq:nabla3-psi-k-bound-1}, 
we obtain that
\begin{align}\label{eq:nabla3-psi-k-bound-2}
\|\bm\nabla^3 \psi_k(\bm{x})\|_1 \lesssim \log_+^{\frac{3}{2}}d  \cdot \frac{e^{-3t}}{\sqrt{(1-e^{-2t})\varepsilon_1}} \cdot \frac{1}{(1-e^{-2t})\varepsilon_1 + e^{-2t}\varepsilon_k},
\end{align}
almost surely.
Finally, by combining \eqref{eq:Lindeberg-swapping}, \eqref{eq:Lindeberg-summand}, \eqref{eq:Lindeberg-summand-3},  \eqref{eq:Lindeberg-summand-4} and \eqref{eq:nabla3-psi-k-bound-2}, and using the tower property of the conditional expectation, we arrive at the bound
\begin{equation}\label{eq:Hd}
\sup_{\mathcal{A} \in \mathscr{H}_d}|\mathbb{E}[h_{\mathcal{A},t}(\bm{S}_n)-h_{\mathcal{A},t}(\bm{T}_1)]| \lesssim \log_+^{\frac{3}{2}}d \frac{e^{-t}}{\sqrt{(1-e^{-2t})\varepsilon_1}}  \sum_{k=1}^n   \mathbb{E} \left[  \frac{ e^{-2t} \|\bm{X}_k\|_{\infty}^3}{(1-e^{-2t})\varepsilon_1 + e^{-2t}\varepsilon_k} \right],
\end{equation}
which by  \eqref{eq:Berry-Esseen-Kol-decompose-2} yields the final master bound
\begin{align}\label{eq:master}
d_{\mathsf{K}}(\bm{S}_n,\bm{T}_1) & \lesssim \log_+^{\frac{3}{2}}d \frac{e^{-t}}{\sqrt{(1-e^{-2t})\varepsilon_1}}  \sum_{k=1}^n   \mathbb{E} \left[  \frac{ e^{-2t} \|\bm{X}_k\|_{\infty}^3}{(1-e^{-2t})\varepsilon_1 + e^{-2t}\varepsilon_k} \right] +\frac{\sqrt{1-e^{-2t}}}{e^{-t}}\sqrt{\frac{d_{\max}(\bm{\Sigma}_n)}{d_{\min}(\bm{\Sigma}_n)}} \log_+d
\end{align}

We now give the proofs of Theorems \ref{thm:Berry-Esseen-Kol-naive} and \ref{thm:Berry-Esseen-Kol} separately. 

\paragraph{Proof of Theorem \ref{thm:Berry-Esseen-Kol-naive}} 
It is sufficient to bound the right hand side of \eqref{eq:Hd}. This can be done straightforwardly using the fact that $e^{-2t} \varepsilon_k \geq 0$ for all $k$, almost surely. Thus, 
\begin{align*}
 \sum_{k=1}^n \mathbb{E}\left[ \frac{ e^{-2t} \|\bm{X}_k\|_{\infty}^3}{(1-e^{-2t})\varepsilon_1 + e^{-2t}\varepsilon_k} \right]&  \leq \frac{e^{-2t}}{(1-e^{-2t}) \varepsilon_1} \sum_{k=1}^n \mathbb{E}\|\bm{X}_k\|_{\infty}^3\\
& = \frac{e^{-2t}}{(1-e^{-2t}) \lambda_{\min}(\bm{\Sigma}_n)} \sum_{k=1}^n \frac{\mathbb{E}\|\bm{X}_k\|_{\infty}^3}{n}.
\end{align*}
In the identity above, we have used the definition of $\varepsilon_1 = n \lambda_{\min}(\bm{\Sigma}_n)$.
Thus, combining the previous bounds with \eqref{eq:Berry-Esseen-Kol-decompose-2},
\[
d_{\mathsf{K}}(\bm{S}_n,\bm{T}_1) \lesssim (\log_+ d)^{\frac{3}{2}}  \frac{1}{\lambda^{3/2}_{\min}(\bm{\Sigma}_n)} \frac{1}{\sqrt{n}}\left( \sum_{k=1}^n \frac{\mathbb{E}\|\bm{X}_k\|_{\infty}^3}{n} \right)  \frac{e^{-3t}}{(1-e^{-2t})^{3/2}} + \frac{\sqrt{1-e^{-2t}}}{e^{-t}} \sqrt{\frac{d_{\max}(\bm{\Sigma}_n)}{d_{\min}(\bm{\Sigma}_n)}} \log_+d.
\]
Theorem \ref{thm:Berry-Esseen-Kol-naive} is recovered by choosing $t > 0$ so that 
\[
\frac{\sqrt{1-e^{-2t}}}{e^{-t}}  = \left( (\log_+ d)^{1/2} \frac{1}{\lambda^{3/2}_{\min}(\bm{\Sigma}_n)} \frac{1}{\sqrt{n}}\left( \sum_{k=1}^n \frac{\mathbb{E}\|\bm{X}_k\|_{\infty}^3}{n} \right) \sqrt{{\frac{d_{\min}(\bm{\Sigma}_n)}{d_{\max}(\bm{\Sigma}_n)}}} \right)^{1/4} .
\]

\paragraph{Proof of Theorem \ref{thm:Berry-Esseen-Kol}} The tighter bound in Theorem \ref{thm:Berry-Esseen-Kol}
 is achieved through a telescoping argument on  the term $\mathbb{E}_{k-1}\|\bm{X}_k\|_{\infty}^3$ in \eqref{eq:Hd}.
Recalling that $\bm{V}_k = \bm{P}_k - \bm{P}_{k+1}$, we have 
\begin{align}\label{eq:3rd-momentum-bound}
\mathbb{E}_{k-1}\|\bm{X}_k\|_{\infty}^3 &\leq M \lambda_{\min}(\bm{P}_k - \bm{P}_{k+1}) \nonumber \\ 
&\leq M [\lambda_{\min}(\bm{P}_k)-\lambda_{\min}(\bm{P}_{k+1})]\nonumber \\ 
&= M (\varepsilon_k - \varepsilon_{k+1})
\end{align}
almost surely, where the first inequality follows from Assumption \eqref{eq:telescope-condition} and the second inequality from the fact that for any symmetric and positive semi-definite matrices $\bm{A}$ and $\bm{B}$ of equal size, $\lambda_{\min}(\bm{A} + \bm{B}) \geq \lambda_{\min}(\bm{A}) + \lambda_{\min}(\bm{B})$. Consequently, each summand on the right-hand-side of \eqref{eq:Hd} can be bounded by
\begin{align*}
&\mathbb{E} \left[ \frac{ e^{-2t}\ \|\bm{X}_k\|_{\infty}^3}{(1-e^{-2t})\varepsilon_1 + e^{-2t}\varepsilon_k} \right] \\ 
&= \mathbb{E}\left[ \frac{e^{-2t}\ \mathbb{E}_{k-1}\|\bm{X}_k\|_{\infty}^3}{(1-e^{-2t})\varepsilon_1 + e^{-2t}\varepsilon_k} \right] \\ 
&\leq M\mathbb{E} \left[ \frac{e^{-2t}(\varepsilon_k - \varepsilon_{k+1})}{(1-e^{-2t})\varepsilon_1 + e^{-2t}\varepsilon_k}\right]
\end{align*}

Since the function $x \in [0,\infty) \mapsto f(x) = \frac{1}{x+(1-e^{-2t})\varepsilon_1}$ is strictly decreasing almost surely, the summation can be bounded by a Riemann integral. Specifically,
\begin{align*}
\sum_{k=1}^n \frac{e^{-2t}(\varepsilon_k - \varepsilon_{k+1})}{(1-e^{-2t})\varepsilon_1 + e^{-2t}\varepsilon_k} &\leq \int_0^{e^{-2t}\varepsilon_1} \frac{1}{x+(1-e^{-2t})\varepsilon_1} \mathrm{d}x \\ 
&= \log \frac{1}{1-e^{-2t}}.
\end{align*}
Plugging into \eqref{eq:Hd}, we obtain
\[
\sup_{\mathcal{A} \in \mathscr{H}_d}|\mathbb{E}[h_{\mathcal{A},t}(\bm{S}_n)-h_{\mathcal{A},t}(\bm{T}_1)]| \lesssim M\log_+^{\frac{3}{2}}d \frac{e^{-t}}{\sqrt{(1-e^{-2t})\varepsilon_1}} \cdot \log \frac{1}{1-e^{-2t}}
\]
almost surely.
Combining with the bound in \eqref{eq:Berry-Esseen-Kol-decompose-2}, we have established that 
\[
d_{\mathsf{K}}(\bm{S}_n,\bm{T}_1) \lesssim M(\log_+ d)^{\frac{3}{2}} \frac{e^{-t}}{\sqrt{(1-e^{-2t})\varepsilon_1}} \log \frac{1}{1-e^{-2t}} + \frac{\sqrt{1-e^{-2t}}}{e^{-t}} \sqrt{\frac{d_{\max}(\bm{\Sigma}_n)}{d_{\min}(\bm{\Sigma}_n)}} \log_+d.
\]

The claim of the theorem 
follows by noting that $\varepsilon_1 = n\lambda_{\min}(\bm{\Sigma}_n)$ and by choosing 
\[
t =  \frac{\log (1 + \Delta^2)}{2},
\]
or, equivalently,
\[
\frac{\sqrt{1-e^{-2t}}}{e^{-t}} = \Delta,
\]


\subsection{Proof of Theorem \ref{thm:Berry-Esseen-Kol-MVB}}\label{app:proof-thm-Berry-Esseen-Kol-MVB}
By Lemma 1 in \cite{arun2020hdCLT}, for any $\varepsilon > 0$,
 \begin{align}\label{eq:decompose-1}
 d_{\mathsf{K}}\left(\frac{1}{\sqrt{n}}\sum_{k=1}^n \bm{X}_k,\mathcal{N}(\bm{0},\bm{\Sigma}_n)\right) &= d_{\mathsf{K}}\left(\sum_{k=1}^n \bm{X}_k ,\mathcal{N}(\bm{0},n\bm{\Sigma}_n )\right)\nonumber\\ 
 &\lesssim d_{\mathsf{K}}\left(\sum_{k=1}^n \bm{X}_k + \varepsilon \bm{Z},\mathcal{N}(\bm{0},n\bm{\Sigma}_n + \varepsilon^2 \bm{I})\right) + \frac{\varepsilon \log_+ d}{\sqrt{nd_{\min}(\bm{\Sigma}_n)}}.
 \end{align}
 Here, $\bm{Z}$ represents a $d$-dimensional standard Gaussian random variable that is independent of $\{\bm{X}_k\}_{k \in [n]}$. We now focus on the first term on the right hand side. For simplicity, we again denote
 \begin{align*}
 \bm{S}_k = \sum_{j=1}^k \bm{X}_j, \quad \bm{Y}_j = \bm{V}_j^{\frac{1}{2}}\bm{Z}_j, \quad \bm{T}_k = \sum_{j=k}^n \bm{Y}_j, \quad \text{and} \quad \bm{T}'_k = \bm{T}_k + \varepsilon \bm{Z},
 \end{align*}
 where $\{\bm{Z}_j\}_{j \in [n]}$ denote a sequence of $i.i.d.$ standard Gaussian random variables in $\mathbb{R}^d$ that are independent of $\{\bm{X}_k\}_{k \in [n]}$ and $\bm{Z}$.
 Furthermore, let 
 \begin{align*}
 \bm{P}_k = \sum_{j=k}^n \bm{V}_j
 \end{align*}
 and recall  that, by assumption, $\bm{P}_1 = n \bm{\Sigma}_n$.

  By the triangle inequality, the first term on the right-hand-side of \eqref{eq:decompose-1} can be further decomposed as
  \begin{align*}
  d_{\mathsf{K}}\left(\sum_{k=1}^n \bm{X}_k + \varepsilon \bm{Z},\mathcal{N}(\bm{0},n\bm{\Sigma}_n + \varepsilon^2 \bm{I})\right) \leq d_{\mathsf{K}}\left(\bm{S}_n + \varepsilon \bm{Z},\bm{T}'_1\right) + d_{\mathsf{K}}\left(\bm{T}'_1,\mathcal{N}(\bm{0},n\bm{\Sigma}_n + \varepsilon^2 \bm{I})\right).
  \end{align*}
  By definition, $\bm{T}'_1 \mid \mathscr{F}_{n} \sim \mathcal{N}(\bm{0},n\bar{\bm{\Sigma}}_n + \varepsilon^2 \bm{I})$, hence by the law of iterated expectations we have that 
  \begin{align*}
  0\leq d_{\mathsf{K}}\left(\bm{T}'_1,\mathcal{N}(\bm{0},n\bm{\Sigma}_n + \varepsilon^2 \bm{I})\right) &= \sup_{\mathcal{A} \in \mathscr{H}_d}\left|\mathbb{P}(\bm{T}'_1 \in \mathcal{A}) - \mathcal{N}\{\bm{0},n\bm{\Sigma}_n+\varepsilon^2 \bm{I}\}(\mathcal{A})\right| \\ 
  &= \sup_{\mathcal{A} \in \mathscr{H}_d}\left|\mathbb{E}[\mathbb{P}(\bm{T}'_1 \in \mathcal{A}\mid \mathscr{F}_n)] - \mathcal{N}\{\bm{0},n\bm{\Sigma}_n+\varepsilon^2 \bm{I}\}(\mathcal{A})\right| \\ 
  &\leq \mathbb{E}\sup_{\mathcal{A} \in \mathscr{H}_d}\left| \mathcal{N}\{\bm{0},n\bar{\bm{\Sigma}}_n+\varepsilon^2 \bm{I}\}(\mathcal{A}) - \mathcal{N}\{\bm{0},n\bm{\Sigma}_n+\varepsilon^2 \bm{I}\}(\mathcal{A})\right| \\ 
  &\leq \mathbb{E} d_{\mathsf{K}}(\mathcal{N}(\bm{0},n\bar{\bm{\Sigma}}_n + \varepsilon^2 \bm{I}),\mathcal{N}(\bm{0},n\bm{\Sigma}_n + \varepsilon^2 \bm{I})) = 0, 
  \end{align*}
  where the last inequality follows from Assumption \ref{as:fixed-variation} that $\bar{\bm{\Sigma}}_n=\bm{\Sigma}_n$.
 As a result, $d_{\mathsf{K}}\left(\bm{T}'_1,\mathcal{N}(\bm{0},n\bm{\Sigma}_n + \varepsilon^2 \bm{I})\right) = 0$. It now boils down to bounding
 \begin{align*}
 d_{\mathsf{K}}\left(\sum_{k=1}^n \bm{X}_k + \varepsilon \bm{Z},\mathcal{N}(\bm{0},n\bm{\Sigma}_n + \varepsilon^2 \bm{I})\right) = d_{\mathsf{K}}\left(\bm{S}_n + \varepsilon \bm{Z},\bm{T}'_1\right).
 \end{align*}
 Towards this end, we invoke Lindeberg's swapping technique:
 \begin{align}\label{eq:th2.2_Lindeberg-swapping}
 d_{\mathsf{K}}\left(\bm{S}_n + \varepsilon \bm{Z},\bm{T}'_1\right) &= \sup_{\mathcal{A} \in \mathscr{H}_d}|\mathbb{P}(\bm{S}_n + \varepsilon \bm{Z} \in \mathcal{A}) - \mathbb{P}(\bm{T}'_1 \in \mathcal{A})| \nonumber\\ 
 &= \sup_{\mathcal{A} \in \mathscr{H}_d} \left|\sum_{k=1}^n \left[\mathbb{P}(\bm{S}_{k-1} + \bm{X}_k + \bm{T}'_{k+1}\in \mathcal{A}) - \mathbb{P}(\bm{S}_{k-1} + \bm{Y}_k + \bm{T}'_{k+1}\in \mathcal{A})\right]\right| \nonumber \\ 
 &\leq \sup_{\mathcal{A} \in \mathscr{H}_d} \sum_{k=1}^n  \left|\mathbb{P}(\bm{S}_{k-1} + \bm{X}_k + \bm{T}'_{k+1}\in \mathcal{A}) - \mathbb{P}(\bm{S}_{k-1} + \bm{Y}_k + \bm{T}'_{k+1}\in \mathcal{A})\right|.
 \end{align}

 For each summand in the last line of \eqref{eq:th2.2_Lindeberg-swapping}, we observe that since $\bm{P}_{k+1}$ is measurable with respect to $\mathscr{F}_{k-1}$, $\bm{T}'_{k+1}\mid \mathscr{F}_{k-1} \sim \mathcal{N}(\bm{0},\bm{P}_{k+1} + \varepsilon^2 \bm{I})$. For simplicity, we denote
 \begin{align*}
 \varepsilon_k^2 := \lambda_{\min}(\bm{P}_{k+1} + \varepsilon^2 \bm{I}),
 \end{align*}
 which is also measurable with respect to $\mathscr{F}_{k-1}$. 
  Meanwhile, condition \eqref{eq:MVB-lower} guarantees that
  \begin{align*}
  \bm{P}_{k+1}  + \varepsilon^2 \bm{I}\succeq ((n-k)\alpha + \varepsilon^2)\bm{I}, \quad  \text{which implies} \quad \varepsilon_k^2 \geq (n-k)\alpha + \epsilon^2.
  \end{align*} 
In this way, we can construct independent standard Gaussian random variables $\bm{\eta},\bm{\eta}'$ such that
 \begin{align*}
 \bm{T}'_{k+1}\mid \mathscr{F}_{k-1} \overset{d}{=} \varepsilon_k\bm{\eta} + (\bm{P}_{k+1} -\lambda_{\min}(\bm{P}_{k+1})\bm{I})^{\frac{1}{2}}\bm{\eta}'.
 \end{align*}
Hence, by the law of iterated expectations, 
 \begin{align}\label{eq:th2.2_Lindeberg-summand}
 &\mathbb{P}(\bm{S}_{k-1} + \bm{X}_k + \bm{T}'_{k+1}\in \mathcal{A}) - \mathbb{P}(\bm{S}_{k-1} + \bm{Y}_k + \bm{T}'_{k+1}\in \mathcal{A})\nonumber \\ 
 &= \mathbb{E}\left[\mathbb{P}(\bm{S}_{k-1} + \bm{X}_k + \bm{T}'_{k+1}\in \mathcal{A}\mid \mathscr{F}_{k-1}) - \mathbb{P}(\bm{S}_{k-1} + \bm{Y}_k + \bm{T}'_{k+1}\in \mathcal{A}\mid \mathscr{F}_{k-1})\right] \nonumber \\ 
 &= \mathbb{E}\Bigg[\mathbb{P}(\bm{S}_{k-1} + \bm{X}_k + \varepsilon_k\bm{\eta} + (\bm{P}_{k+1} -\lambda_{\min}(\bm{P}_{k+1})\bm{I})^{\frac{1}{2}}\bm{\eta}'\in \mathcal{A}\mid \mathscr{F}_{k-1})\nonumber \\ 
 &\quad - \mathbb{P}(\bm{S}_{k-1} + \bm{Y}_k + \varepsilon_k\bm{\eta} + (\bm{P}_{k+1} -\lambda_{\min}(\bm{P}_{k+1})\bm{I})^{\frac{1}{2}}\bm{\eta}'\in \mathcal{A}\mid \mathscr{F}_{k-1})\Bigg] \nonumber \\ 
 &= \mathbb{E}\Bigg[\mathbb{P}(\bm{S}_{k-1} + \bm{X}_k + \varepsilon_k\bm{\eta} + (\bm{P}_{k+1} -\lambda_{\min}(\bm{P}_{k+1})\bm{I})^{\frac{1}{2}}\bm{\eta}'\in \mathcal{A}\mid \mathscr{F}_{k-1}\vee \sigma(\bm{\eta}')) \nonumber \\ 
 &\quad - \mathbb{P}(\bm{S}_{k-1} + \bm{Y}_k + \varepsilon_k\bm{\eta} + (\bm{P}_{k+1} -\lambda_{\min}(\bm{P}_{k+1})\bm{I})^{\frac{1}{2}}\bm{\eta}'\in \mathcal{A}\mid \mathscr{F}_{k-1}\vee \sigma(\bm{\eta}'))\Bigg].
 \end{align}
 For simplicity, we define 
 \begin{align*}
 \bm{S}'_{k} = \bm{S}_k +  (\bm{P}_{k+1} -\lambda_{\min}(\bm{P}_{k+1})\bm{I})^{\frac{1}{2}}\bm{\eta}'.
 \end{align*}
In this way, \eqref{eq:th2.2_Lindeberg-summand} simplifies to
 \begin{align}\label{eq:Lindeberg-summand-simplify}
 &\mathbb{P}(\bm{S}_{k-1} + \bm{X}_k + \bm{T}'_{k+1}\in \mathcal{A}) - \mathbb{P}(\bm{S}_{k-1} + \bm{Y}_k + \bm{T}'_{k+1}\in \mathcal{A}) \nonumber \\ 
 &= \mathbb{E}[\varphi_{\mathcal{A},\varepsilon_k}(\bm{S}'_{k-1} + \bm{X}_k) - \varphi_{\mathcal{A},\varepsilon_k}(\bm{S}'_{k-1} + \bm{Y}_k) \mid \mathscr{F}_{k-1}\vee \sigma(\bm{\eta}'))]
 \end{align}
For every $1 \leq k < n$, by Taylor's expansion 
 \begin{align*}
 \varphi_{\mathcal{A},\varepsilon_k}(\bm{S}'_{k-1} + \bm{X}_k) &= \varphi_{\mathcal{A},\varepsilon_k}(\bm{S}'_{k-1}) + \nabla \varphi_{\mathcal{A},\varepsilon_k}(\bm{S}'_{k-1}) [\bm{X}_k] + \frac{1}{2}\nabla^2 \varphi_{\mathcal{A},\varepsilon_k}(\bm{S}'_{k-1}) [\bm{X}_k]^2 \\ 
 &+  \frac{1}{6}\int_{\theta=0}^1 \langle \nabla^3 \varphi_{\mathcal{A},\varepsilon_k}(\bm{S}'_{k-1} + \theta \bm{X}_k), \bm{X}_k^{\otimes 3}\rangle, \quad \text{and} \\ 
 \varphi_{\mathcal{A},\varepsilon_k}(\bm{S}'_{k-1} + \bm{Y}_k) &= \varphi_{\mathcal{A},\varepsilon_k}(\bm{S}'_{k-1}) + \nabla \varphi_{\mathcal{A},\varepsilon_k}(\bm{S}'_{k-1}) [\bm{Y}_k] + \frac{1}{2}\nabla^2 \varphi_{\mathcal{A},\varepsilon_k}(\bm{S}'_{k-1}) [\bm{Y}_k]^2 \\ 
 &+  \frac{1}{6}\int_{\theta=0}^1 \langle \nabla^3 \varphi_{\mathcal{A},\varepsilon_k}(\bm{S}'_{k-1} + \theta \bm{Y}_k), \bm{Y}_k^{\otimes 3}\rangle.
 \end{align*}
 Here, since $\bm{S}'_{k-1}$ is measurable with respect to the $\sigma$-field $\mathscr{F}_{k-1}\vee \sigma(\bm{\eta}')$, and because $\mathbb{E}[\bm{X}_k\bm{X}_k^\top\mid \mathscr{F}_{k-1}] = \mathbb{E}[\bm{Y}_k\bm{Y}_k^\top\mid \mathscr{F}_{k-1}] = \bm{V}_k$, the difference in \eqref{eq:Lindeberg-summand-simplify} can be bounded by 
 \begin{align}\label{eq:Lindeberg-summand-bound-1}
 &\left|\mathbb{E}[\varphi_{\mathcal{A},\varepsilon_k}(\bm{S}'_{k-1} + \bm{X}_k) - \varphi_{\mathcal{A},\varepsilon_k}(\bm{S}'_{k-1} + \bm{Y}_k) \mid \mathscr{F}_{k-1}\vee \sigma(\bm{\eta}'))]\right| \nonumber \\ 
 &\leq \sup_{x \in \mathbb{R}^d} \|\nabla^3 \varphi_{\mathcal{A},\varepsilon_k}(\bm{x})\|_1 \cdot \mathbb{E}[\|\bm{X}_k\|_{\infty}^3 + \|\bm{Y}_k\|_{\infty}^3 \mid \mathscr{F}_{k-1}];
 \end{align}
 In this expression, Theorem \ref{thm:Fang} guarantees that $\sup_{x \in \mathbb{R}^d} \|\nabla^3 \varphi_{\mathcal{A},\sigma}(\bm{x})\|_1 \lesssim \varepsilon_k^{-3}\log_+^{\frac{3}{2}}d$; therefore, the first $n-1$ summands in \eqref{eq:th2.2_Lindeberg-swapping} are bounded by
 \begin{align*}
 &\left|\mathbb{P}(\bm{S}_{k-1} + \bm{X}_k + \bm{T}'_{k+1}\in \mathcal{A}) - \mathbb{P}(\bm{S}_{k-1} + \bm{Y}_k + \bm{T}'_{k+1}\in \mathcal{A})\right| \\ 
 &\leq \varepsilon_k^{-3}\log_+^{\frac{3}{2}}d \cdot \mathbb{E}[\|\bm{X}_k\|_{\infty}^3 + \|\bm{Y}_k\|_{\infty}^3 \mid \mathscr{F}_{k-1}];
 \end{align*}
 For the last summand ($k=n$), we simply invoke the first-order expansion to obtain
 \begin{align*}
 &\left|\mathbb{E}[\varphi_{\mathcal{A},\varepsilon_n}(\bm{S}'_{n-1} + \bm{X}_n) - \varphi_{\mathcal{A},\varepsilon_n}(\bm{S}'_{n-1} + \bm{Y}_n) \mid \mathscr{F}_{n-1}\vee \sigma(\bm{\eta}'))]\right| \\ 
 &\leq \sup_{x \in \mathbb{R}^d} \|\nabla \varphi_{\mathcal{A},\varepsilon_n}(\bm{x})\|_1 \cdot \mathbb{E}[\|\bm{X}_n\|_{\infty} + \|\bm{Y}_n\|_{\infty} \mid \mathscr{F}_{n-1}] \\ 
 &\leq \varepsilon^{-1} \log_+^{\frac{1}{2}}d\cdot \mathbb{E}[\|\bm{X}_n\|_{\infty} + \|\bm{Y}_n\|_{\infty} \mid \mathscr{F}_{n-1}],
 \end{align*}
 where the last line follows from the fact that $\varepsilon_n = \varepsilon$ and Theorem \ref{thm:Fang}.

 In combination, we have the following proposition:
 \begin{proposition}\label{prop:high-dim-Rollin}
 It can be guaranteed that for every $\varepsilon > 0$,
 \begin{align}
 d_{\mathsf{K}}\left(\bm{S}_n + \varepsilon \bm{Z},\bm{T}'_1\right) 
 &\leq \sum_{k=1}^{n-1} \mathbb{E}\left[\varepsilon_k^{-3}\log_+^{\frac{3}{2}}d \cdot \mathbb{E}[\|\bm{X}_k\|_{\infty}^3 + \|\bm{Y}_k\|_{\infty}^3 \mid \mathscr{F}_{k-1}]\right] \nonumber \\ 
 &+ \varepsilon^{-1} \log_+^{\frac{1}{2}}d\cdot\mathbb{E}[\|\bm{X}_n\|_{\infty} + \|\bm{Y}_n\|_{\infty} \mid \mathscr{F}_{n-1}] 
 \end{align}
 \end{proposition}

 Next, we discuss the bound on the Kolmogorov distance under (i) only the condition \eqref{eq:MVB-upper}, and (ii) both condition \eqref{eq:MVB-upper} and condition \eqref{eq:MVB-lower}.
 
 \paragraph{If only condition \eqref{eq:MVB-upper} is guaranteed:} Simply observe $\varepsilon_k \geq \varepsilon$ for all $k \in [n]$ almost surely, and that Theorem \ref{thm:Gaussian-inf-norm} implies
 \begin{align*}
\mathbb{E}\|\bm{Y}_k\|_{\infty}^3 \lesssim \beta^{\frac{3}{2}}\log_+^{\frac{3}{2}}d.
 \end{align*}
 
Therefore, Proposition \ref{prop:high-dim-Rollin} and the tower property of conditional expectations directly implies
 \begin{align*}
d_{\mathsf{K}}\left(\bm{S}_n + \varepsilon \bm{Z},\bm{T}'_1\right) &\leq \varepsilon^{-3}\log_+^{\frac{3}{2}} d \sum_{k=1}^n \mathbb{E}\left[\mathbb{E}[\|\bm{X}_k\|_{\infty}^3 + \|\bm{Y}_k\|_{\infty}^3 \mid \mathscr{F}_{k-1}]\right] \\ 
& = \varepsilon^{-3} \log_+^{\frac{3}{2}} d\sum_{k=1}^n \mathbb{E}[\|\bm{X}_k\|_{\infty}^3 + \|\bm{Y}_k\|_{\infty}^3] \\ 
&\lesssim \varepsilon^{-3} \cdot n(\gamma+\beta)^{\frac{3}{2}}\log_+^{3} d;
 \end{align*}
 combining this bound with \eqref{eq:decompose-1}, we obtain
\begin{align*}
 d_{\mathsf{K}}\left(\frac{1}{\sqrt{n}}\sum_{k=1}^n \bm{X}_k,\mathcal{N}(\bm{0},\bm{\Sigma}_n)\right) \lesssim \varepsilon^{-3} \cdot n(\gamma+\beta)^{\frac{3}{2}}\log_+^3 d + \frac{\varepsilon \log_+ d}{\sqrt{nd_{\min}(\bm{\Sigma}_n)}};
\end{align*}
the upper bound \eqref{eq:Berry-Esseen-Kol-MVB1} follows by taking 
\begin{align*}
\varepsilon = n^{\frac{3}{8}} (\gamma+\beta)^{\frac{3}{8}} d_{\min}^{\frac{1}{8}}(\bm{\Sigma}_n) \log_+^{\frac{1}{2}} d.
\end{align*}

 \paragraph{If both conditions \eqref{eq:MVB-upper} and \eqref{eq:MVB-lower} are satisfied:} We have that, almost surely, 
 \begin{align*}
 \varepsilon_k^2 \geq (n-k)\alpha + \varepsilon^2.
 \end{align*}
 Therefore, proposition \ref{prop:high-dim-Rollin} can be simplified, again using the tower property, as 
 \begin{align*}
 d_{\mathsf{K}}\left(\bm{S}_n + \varepsilon \bm{Z},\bm{T}'_1\right) 
 &\leq\sum_{k=1}^{n-1} ((n-k)\alpha + \varepsilon^2)^{-\frac{3}{2}}\log_+^{\frac{3}{2}}d \cdot (\gamma^{\frac{3}{2}} + \beta^{\frac{3}{2}})\log_+^{\frac{3}{2}}d \\ 
 &+ \varepsilon^{-1} \log_+^{\frac{1}{2}}d\cdot (\gamma^{\frac{1}{2}} + \beta^{\frac{1}{2}})\log_+d.
 \end{align*}
 Here, the summation can be bounded by an integral as
 \begin{align*}
 \sum_{k=1}^{n-1} ((n-k)\alpha + \varepsilon^2)^{-\frac{3}{2}} \leq \frac{1}{\alpha}\int_{\varepsilon^2}^{\infty} x^{-\frac{3}{2}}\mathrm{d}x = \frac{2}{\alpha \varepsilon}.
 \end{align*}
 Hence we obtain
 \begin{align*}
 d_{\mathsf{K}}\left(\bm{S}_n + \varepsilon \bm{Z},\bm{T}'_1\right) & \lesssim \frac{2(\gamma^{\frac{3}{2}} + \beta^{\frac{3}{2}})}{\alpha \varepsilon} \log_+^3 d + \frac{\gamma^{\frac{1}{2}} + \beta^{\frac{1}{2}}}{\varepsilon} \log_+ d \lesssim \frac{(\gamma+\beta)^{\frac{3}{2}}}{\alpha \varepsilon} \log_+^3 d .
 \end{align*}
 In combination with \eqref{eq:decompose-1} we obtain that
 \begin{align*}
 d_{\mathsf{K}}\left(\frac{1}{\sqrt{n}}\sum_{k=1}^n \bm{X}_k,\mathcal{N}(\bm{0},\bm{\Sigma}_n)\right) \lesssim \frac{(\gamma+\beta)^{\frac{3}{2}}}{\alpha \varepsilon} \log_+^3 d + \frac{\varepsilon \log_+ d}{\sqrt{nd_{\min}(\bm{\Sigma}_n)}},
 \end{align*}
 for any $\varepsilon>0$. By taking
 \begin{align*}
 \varepsilon = \sqrt{\frac{(\gamma+\beta)^{\frac{3}{2}}}{\alpha d_{\min}^{\frac{1}{2}}(\bm{\Sigma}_n)}} \log_+ d n^{\frac{1}{4}},
 \end{align*}
 the theorem follows immediately.

\subsection{Proof of Theorem \ref{thm:Berry-Esseen-Kol-Markov}}
This proof is inspired by the multi-dimensional Yurinskii's coupling technique invoked by \cite{belloni2018highdimensionalcentrallimit,yurinksi.martingale} and \cite{wu2025uncertainty}. Borrowing the notation of Theorem \ref{thm:Berry-Esseen-Kol}, we define
\begin{align*}
&\bm{V}_k = \mathbb{E}[\bm{f}(s_{k-1},s_k)\bm{f}^\top(s_{k-1},s_k)\mid \mathscr{F}_{k-1}] = \mathbb{E}[\bm{f}(s_{k-1},s_k)\bm{f}^\top(s_{k-1},s_k)\mid s_{k-1}], \\ 
&\bar{\bm{\Sigma}}_n = \frac{1}{n}\sum_{k=1}^n \bm{V}_k,  \quad \text{and} \quad \bm{S}_k = \sum_{j=1}^n \bm{f}(s_{j-1},s_j).
\end{align*}
The proof consists of the following steps:
\begin{enumerate}
    \item Define a quantity $\kappa \asymp \log(nd) n^{-\frac{1}{2}}$, such that with high probability,
    \begin{align*}
        \|\bar{\bm{\Sigma}}_n - \bm{\Sigma}\| \leq \kappa.
    \end{align*}
    \item Construct a martingale $\{\tilde{\bm{S}}_k\}_{k=1}^{n+1}$, such that
    \begin{align*}
        \mathbb{E}[\tilde{\bm{S}}_{n+1}\tilde{\bm{S}}_{n+1}^\top] = n(\bm{\Sigma} + \kappa \bm{I}) \quad \text{a.s.,}
    \end{align*}
    and that $\tilde{\bm{S}}_{n+1}$ is close to $\bm{S}_n$.
    \item Apply a generalized version of Theorem \ref{thm:Berry-Esseen-Kol} on $\tilde{\bm{S}}_{n+1}$;
    \item Combine the terms and complete the proof by invoking anti-concentration and Gaussian comparison results with regard to the Kolmogorov distance.
\end{enumerate}
\paragraph{Step 1: Specify $\kappa$.} By definition, the condition that $\alpha \bm{I} \preceq \mathbb{E}_{s' \sim P(\cdot \mid s)}[\bm{f}(s,s')\bm{f}^\top(s,s')] \preceq \beta\bm{I}$ for all $s \in \mathcal{S}$ effectively guarantees that
\begin{align*}
\|\mathbb{E}_{s' \sim P(\cdot \mid s)}[\bm{f}(s,s')\bm{f}^\top(s,s')] - \bm{\Sigma}\| \leq \beta - \alpha, \quad \forall s \in \mathcal{S}.
\end{align*}
As a direct consequence, the matrix Hoeffding's inequality on Markov chains \cite[Theorem 2.5]{neeman2024concentration}, \cite[Corollary 1]{wu2025uncertainty} indicates
\begin{align}\label{eq:defn-kappa}
\|\bar{\bm{\Sigma}}_n - \bm{\Sigma}\| \lesssim \sqrt{\frac{q}{1-\lambda} \log \left(2d n \left\|\frac{\mathrm{d}\nu}{\mathrm{d}\mu}\right\|_{\mu,p}\right)} \cdot \frac{\beta-\alpha}{\sqrt{n}}=:\kappa
\end{align}
with probability at least $1-n^{-\frac{1}{2}}$.
\paragraph{Step 2: Construct the martingale $\{\tilde{\bm{S}}_k\}_{k \in [n+1]}$.} Define the stopping time
\begin{align*}
\tau := \sup\left\{t \leq n: \sum_{i=1}^t \bm{V}_i \preceq n(\bm{\Sigma}+\kappa \bm{I})\right\},
\end{align*}
and define a sequence $\{\bm{X}_j\}_{j \in [n+1]}$ as
\begin{align*}
\bm{X}_j = \begin{cases}
\bm{f}(s_{j-1},s_j), \quad &\text{if} \quad 1 \leq j \leq \tau;\\
\bm{0}, \quad &\text{if} \quad \tau+1 \leq j \leq n; \\ 
(n(\bm{\Sigma}+\kappa \bm{I}) - \sum_{i=1}^{\tau} \bm{V}_i)^{\frac{1}{2}}\bm{\eta}, \quad &\text{if} \quad j= n+1.
\end{cases}
\end{align*}
Here, $\bm{\eta}$ is a $d$-dimensional standard Gaussian random variable independent of $\mathscr{F}_n$. 
By definition, it is easy to verify that  $\{\bm{X}_j\}_{j \in [n+1]}$ is a martingale difference process, and therefore its partial sum
\begin{align*}
\tilde{\bm{S}}_k = \sum_{j=1}^k \bm{X}_j, \quad \forall k \in [n+1]
\end{align*}
is a martingale process.
We now aim to bound the difference $\tilde{\bm{S}}_{n+1}$ and $\bm{S}_n$. Specifically, for any $x \geq 0$, 
\begin{align}\label{eq:Yurinskii-compare-1}
&\mathbb{P}\left(\left\|\frac{\bm{S}_n}{\sqrt{n}}- \frac{\tilde{\bm{S}}_{n+1}}{\sqrt{n}}\right\|_{\infty} \geq x\right) \nonumber \\ 
&\leq \mathbb{P}\left(\left\|\frac{\bm{S}_n}{\sqrt{n}}- \frac{\tilde{\bm{S}}_{n+1}}{\sqrt{n}}\right\|_{\infty} \geq x, \|\bar{\bm{\Sigma}}_n - \bm{\Sigma}\| \leq \kappa \right) + \mathbb{P}( \|\bar{\bm{\Sigma}}_n - \bm{\Sigma}\| \geq \kappa) \nonumber \\ 
&\leq  \mathbb{P}\left(\left\|\frac{\bm{S}_n}{\sqrt{n}}- \frac{\tilde{\bm{S}}_{n+1}}{\sqrt{n}}\right\|_{\infty} \geq x \bigg | \|\bar{\bm{\Sigma}}_n - \bm{\Sigma}\| \leq \kappa \right)  \mathbb{P}( \|\bar{\bm{\Sigma}}_n - \bm{\Sigma}\| \leq \kappa) + \mathbb{P}( \|\bar{\bm{\Sigma}}_n - \bm{\Sigma}\| \geq \kappa)\nonumber \\ 
&\leq \mathbb{P}\left(\left\|\frac{\bm{S}_n}{\sqrt{n}}- \frac{\tilde{\bm{S}}_{n+1}}{\sqrt{n}}\right\|_{\infty} \geq x \bigg | \|\bar{\bm{\Sigma}}_n - \bm{\Sigma}\| \leq \kappa \right) + \mathbb{P}( \|\bar{\bm{\Sigma}}_n - \bm{\Sigma}\| \geq \kappa)
\end{align}
where the second term on the last line is bounded by $n^{-\frac{1}{2}}$ according to \eqref{eq:defn-kappa}. For the first term, it follows by definition that when $\|\bar{\bm{\Sigma}}_n - \bm{\Sigma}\| \leq \kappa$, 
\begin{align*}
\sum_{i=1}^n \bm{V}_i = n\bar{\bm{\Sigma}}_n \preceq n(\bm{\Sigma}_n + \kappa \bm{I}),
\end{align*}
and therefore $\tau = n$. Consequently, the difference between $\bm{S}_n - \tilde{\bm{S}}_{n+1}$ is a Gaussian random variable whose variance matrix 
\begin{align*}
n(\bm{\Sigma} + \kappa \bm{I} - \bar{\bm{\Sigma}}_n) \preceq 2n\kappa \bm{I}.
\end{align*}
Therefore, a direct application of Theorem \ref{thm:Gaussian-inf-norm} indicates
\begin{align*}
\mathbb{P}\left(\left\|\frac{\bm{S}_n}{\sqrt{n}}- \frac{\tilde{\bm{S}}_{n+1}}{\sqrt{n}}\right\|_{\infty} \gtrsim \sqrt{\kappa \log(nd)} \Bigg | \|\bar{\bm{\Sigma}}_n - \bm{\Sigma}\| \leq \kappa \right) \leq n^{-\frac{1}{2}}.
\end{align*}
Hence by taking $x \asymp \sqrt{\kappa \log(nd)}$ in \eqref{eq:Yurinskii-compare-1}, we obtain
\begin{align}\label{eq:Yurinskii-compare-2}
\mathbb{P}\left(\left\|\frac{\bm{S}_n}{\sqrt{n}}- \frac{\tilde{\bm{S}}_{n+1}}{\sqrt{n}}\right\|_{\infty} \gtrsim \sqrt{\kappa \log(nd)}\right) \leq 2n^{-\frac{1}{2}}.
\end{align}
\paragraph{Step 3: Apply generalized Theorem \ref{thm:Berry-Esseen-Kol} on $\tilde{\bm{S}}_{n+1}$.}


Although $\tilde{\bm{S}}_{n+1}$ has a fixed terminal quadratic variation, the condition \eqref{eq:MVB-lower} is not always true when $\tau < n$. Therefore, Theorem \ref{thm:Berry-Esseen-Kol} cannot be directly applied, and a refined analysis is required. Similar to the proof of Theorem \ref{thm:Berry-Esseen-Kol}, we denote, for every $j,k \in [n+1]$ and $\varepsilon>0$,
\begin{align*}
\tilde{\bm{V}}_j = \mathbb{E}[\tilde{\bm{X}}_j \tilde{\bm{X}}_j^\top], \quad \tilde{\bm{Y}}_j = \tilde{\bm{V}}_j^{\frac{1}{2}}\bm{Z}_j,\quad \tilde{\bm{T}}_k = \sum_{j=k}^{n+1}\tilde{\bm{Y}}_j, \quad \text{and} \quad \tilde{\bm{T}}_k' = \bm{T}_k + \varepsilon \bm{Z}.
\end{align*}
In this way, the Kolmogorov distance between $\frac{\tilde{\bm{S}}_{n+1}}{\sqrt{n}}$ and $\mathcal{N}(\bm{0},\bm{\Sigma} + \kappa \bm{I})$ is bounded by 
\begin{align}\label{eq:Markov-auxiliary-Lindeberg}
&d_{\mathsf{K}}\left(\frac{\tilde{\bm{S}}_{n+1}}{\sqrt{n}},\mathcal{N}(\bm{0},\bm{\Sigma} + \kappa \bm{I})\right) \nonumber \\ 
&\lesssim \sup_{\mathcal{A} \in \mathscr{H}_d} \sum_{k=1}^{n+1} |\mathbb{P}(\tilde{\bm{S}}_{k-1} + \tilde{\bm{X}}_k+\tilde{\bm{T}}_{k+1}' \in \mathcal{A}) - \mathbb{P}(\tilde{\bm{S}}_{k-1} + \tilde{\bm{Y}}_k+\tilde{\bm{T}}_{k+1}' \in \mathcal{A})| + \frac{\varepsilon \log_+ d}{\sqrt{n}},
\end{align}
where we invoked the fact that the diagonal elements of $\bm{\Sigma}$ are all equal to 1. Notice that by definition, when $\tau < k < n+1$, we have $\tilde{\bm{X}}_k \equiv \tilde{\bm{Y}}_k \equiv \bm{0}$; and when $k = n+1$, we have $\tilde{\bm{X}}_k \overset{d}{=}\tilde{\bm{Y}}_k \mid \mathscr{F}_{n}$ since they are both Gaussian random variables with mean $\bm{0}$ and variance
\begin{align*}
n(\bm{\Sigma} + \kappa\bm{I})-\sum_{i=1}^{\tau}\bm{V}_i.
\end{align*}
Therefore, whenever $k > \tau$, it can be guaranteed that
\begin{align*}
&\mathbb{P}(\tilde{\bm{S}}_{k-1} + \tilde{\bm{X}}_k+\tilde{\bm{T}}_{k+1}' \in \mathcal{A}) - \mathbb{P}(\tilde{\bm{S}}_{k-1} + \tilde{\bm{Y}}_k+\tilde{\bm{T}}_{k+1}' \in \mathcal{A}) \\ 
&= \mathbb{E}\mathbb{P}(\tilde{\bm{S}}_{k-1} + \tilde{\bm{X}}_k+\tilde{\bm{T}}_{k+1}' \in \mathcal{A}\mid\mathscr{F}_{k-1}) -\mathbb{E}\mathbb{P}(\tilde{\bm{S}}_{k-1} + \tilde{\bm{Y}}_k+\tilde{\bm{T}}_{k+1}' \in \mathcal{A}\mid\mathscr{F}_{k-1}) =0.
\end{align*}
As a direct consequence, the summation on the second line of \eqref{eq:Markov-auxiliary-Lindeberg} can be specified as
\begin{align*}
&\sum_{k=1}^{n+1} |\mathbb{P}(\tilde{\bm{S}}_{k-1} + \tilde{\bm{X}}_k+\tilde{\bm{T}}_{k+1}' \in \mathcal{A}) - \mathbb{P}(\tilde{\bm{S}}_{k-1} + \tilde{\bm{Y}}_k+\tilde{\bm{T}}_{k+1}' \in \mathcal{A})| \\ 
&= \sum_{k=1}^{\tau} |\mathbb{P}(\tilde{\bm{S}}_{k-1} + \tilde{\bm{X}}_k+\tilde{\bm{T}}_{k+1}' \in \mathcal{A}) - \mathbb{P}(\tilde{\bm{S}}_{k-1} + \tilde{\bm{Y}}_k+\tilde{\bm{T}}_{k+1}' \in \mathcal{A})|.
\end{align*}
For every $k \in [\tau]$, define $\tilde{\varepsilon}_k^2 = \varepsilon^2 + (\tau-k)\alpha$. Similar to the proof of Theorem \ref{thm:Berry-Esseen-Kol}, it is easy to verify that
\begin{align*}
|\mathbb{P}(\tilde{\bm{S}}_{k-1} + \tilde{\bm{X}}_k+\tilde{\bm{T}}_{k+1}' \in \mathcal{A}) - \mathbb{P}(\tilde{\bm{S}}_{k-1} + \tilde{\bm{Y}}_k+\tilde{\bm{T}}_{k+1}' \in \mathcal{A})| \lesssim \tilde{\varepsilon}_k^{-3}\log_+^3d(\gamma^{\frac{3}{2}} + \beta^{\frac{3}{2}})
\end{align*}
when $k < \tau$, and
\begin{align*}
|\mathbb{P}(\tilde{\bm{S}}_{\tau-1} + \tilde{\bm{X}}_\tau+\tilde{\bm{T}}_{\tau+1}' \in \mathcal{A}) - \mathbb{P}(\tilde{\bm{S}}_{\tau-1} + \tilde{\bm{Y}}_\tau+\tilde{\bm{T}}_{\tau+1}' \in \mathcal{A})| \lesssim \frac{\gamma^{\frac{1}{2}}+\beta^{\frac{1}{2}}}{\varepsilon}\log^{\frac{1}{2}}_+d.
\end{align*}
Since for any $\tau \leq n$, it can be guaranteed that
\begin{align*}
\sum_{k=1}^{\tau-1}\tilde{\varepsilon}_k^{-3} \leq \frac{2(\gamma^{\frac{3}{2}} + \beta^{\frac{3}{2}})}{\alpha},
\end{align*}
a similar argument to Theorem \ref{thm:Berry-Esseen-Kol} yields
\begin{align}\label{eq:auxiliary-Berry-Esseen}
d_{\mathsf{K}}\left(\frac{\tilde{\bm{S}}_{n+1}}{\sqrt{n}},\mathcal{N}(\bm{0}, \bm{\Sigma} + \kappa \bm{I})\right) \lesssim \frac{\gamma^{\frac{3}{4}} + \beta^{\frac{3}{4}}}{\sqrt{\alpha}}\log_+^2d n^{-\frac{1}{4}}.
\end{align}

\paragraph{Step 4: Complete the proof.} By the triangle inequality,
\begin{align}\label{eq:Markov-Berry-Esseen-decompose-1}
d_{\mathsf{K}}\left(\frac{\bm{S}_n}{\sqrt{n}},\mathcal{N}(\bm{0},\bm{\Sigma})\right) \leq d_{\mathsf{K}}\left(\frac{\bm{S}_n}{\sqrt{n}},\mathcal{N}(\bm{0},\bm{\Sigma} + \kappa \bm{I})\right) + d_{\mathsf{K}}\left(\mathcal{N}(\bm{0},\bm{\Sigma}),\mathcal{N}(\bm{0},\bm{\Sigma} + \kappa \bm{I})\right).
\end{align}
For the second term on the right hand side, Theorem \ref{thm:Gaussian-comparison-Kol} directly implies
\begin{align}\label{eq:Markov-Berry-Esseen-Gaussian-compare}
d_{\mathsf{K}}\left(\mathcal{N}(\bm{0},\bm{\Sigma}),\mathcal{N}(\bm{0},\bm{\Sigma} + \kappa \bm{I})\right) \lesssim \|\bm{\Sigma}^{-1}\| \kappa \log \frac{1}{\kappa }\log_+d.
\end{align}
Notice that this depends on our assumption that the diagonal elements of $\bm{\Sigma}$ are equal to 1. For the first term, we observe for any $\mathcal{A} \in \mathscr{H}_d$ and $x > 0$,
\begin{align}\label{eq:Wu-decompose-upper}
\mathbb{P}\left(\frac{\bm{S}_n}{\sqrt{n}} \in \mathcal{A}\right)&= \mathbb{P}\left(\frac{\bm{S}_n}{\sqrt{n}} \in \mathcal{A}, \left\|\frac{\bm{S}_n - \tilde{\bm{S}}_{n+1}}{\sqrt{n}}\right\|_\infty \leq x\right) + \mathbb{P}\left(\frac{\bm{S}_n}{\sqrt{n}} \in \mathcal{A}, \left\|\frac{\bm{S}_n - \tilde{\bm{S}}_{n+1}}{\sqrt{n}}\right\|_\infty > x\right) \nonumber \\ 
&\leq \mathbb{P}\left(\frac{\tilde{\bm{S}}_{n+1}}{\sqrt{n}} \in \mathcal{A}^x \right) + \mathbb{P}\left(\left\|\frac{\bm{S}_n - \tilde{\bm{S}}_{n+1}}{\sqrt{n}}\right\|_\infty > x\right)\nonumber \\
&\leq \mathcal{N}(\bm{0},\bm{\Sigma} + \kappa \bm{I})\{\mathcal{A}^x\} + d_{\mathsf{K}}\left(\frac{\tilde{\bm{S}}_{n+1}}{\sqrt{n}},\mathcal{N}(\bm{0},\bm{\Sigma} + \kappa \bm{I})\right)+ \mathbb{P}\left(\left\|\frac{\bm{S}_n - \tilde{\bm{S}}_{n+1}}{\sqrt{n}}\right\|_\infty > x\right) \nonumber \\ 
&\leq \mathcal{N}(\bm{0},\bm{\Sigma} + \kappa \bm{I})\{\mathcal{A}\} + \frac{\log_+d}{1+\kappa} x + d_{\mathsf{K}}\left(\frac{\tilde{\bm{S}}_{n+1}}{\sqrt{n}},\mathcal{N}(\bm{0},\bm{\Sigma} + \kappa \bm{I})\right)+ \mathbb{P}\left(\left\|\frac{\bm{S}_n - \tilde{\bm{S}}_{n+1}}{\sqrt{n}}\right\|_\infty > x\right).
\end{align}
Notice that in the last line, we invoked Theorem \ref{thm:anti-concentration} and the fact that all diagonal elements in $\bm{\Sigma}$ are equal to $1$. Hence, by taking $x \asymp \sqrt{\kappa \log(nd)}$ and combining with \eqref{eq:Yurinskii-compare-2} and \eqref{eq:auxiliary-Berry-Esseen}, we obtain
\begin{align}\label{eq:Markov-upper-bound}
\mathbb{P}\left(\frac{\bm{S}_n}{\sqrt{n}} \in \mathcal{A}\right) - \mathcal{N}(\bm{0},\bm{\Sigma} + \kappa \bm{I})\{\mathcal{A}\} 
&\lesssim \sqrt{\kappa \log(nd)} \log_+d + \frac{\gamma^{\frac{3}{4}} + \beta^{\frac{3}{4}}}{\sqrt{\alpha}}\log_+^2d n^{-\frac{1}{4}}. 
\end{align}
Meanwhile, a union bound argument yields
\begin{align*}
\mathbb{P}\left(\frac{\bm{S}_n}{\sqrt{n}} \in \mathcal{A}\right) &\geq  \mathbb{P}\left(\frac{\bm{S}_n}{\sqrt{n}} \in \mathcal{A}, \left\|\frac{\bm{S}_n - \tilde{\bm{S}}_{n+1}}{\sqrt{n}}\right\|_\infty > x\right) -  \mathbb{P}\left( \left\|\frac{\bm{S}_n - \tilde{\bm{S}}_{n+1}}{\sqrt{n}}\right\|_\infty > x\right) \\ 
&\geq \mathbb{P}\left(\frac{\tilde{\bm{S}}_{n+1}}{\sqrt{n}} \in \mathcal{A}^{-x} \right)-  \mathbb{P}\left( \left\|\frac{\bm{S}_n - \tilde{\bm{S}}_{n+1}}{\sqrt{n}}\right\|_\infty > x\right) \\ 
&\geq \mathcal{N}(\bm{0},\bm{\Sigma} + \kappa \bm{I})\{\mathcal{A}^{-x}\} - d_{\mathsf{K}}\left(\frac{\tilde{\bm{S}}_{n+1}}{\sqrt{n}},\mathcal{N}(\bm{0},\bm{\Sigma} + \kappa \bm{I})\right)-  \mathbb{P}\left( \left\|\frac{\bm{S}_n - \tilde{\bm{S}}_{n+1}}{\sqrt{n}}\right\|_\infty > x\right) \\ 
&\geq \mathcal{N}(\bm{0},\bm{\Sigma} + \kappa \bm{I})\{\mathcal{A}\} - \frac{\log_+d}{1+\kappa} x - d_{\mathsf{K}}\left(\frac{\tilde{\bm{S}}_{n+1}}{\sqrt{n}},\mathcal{N}(\bm{0},\bm{\Sigma} + \kappa \bm{I})\right)- \mathbb{P}\left(\left\|\frac{\bm{S}_n - \tilde{\bm{S}}_{n+1}}{\sqrt{n}}\right\|_\infty > x\right).
\end{align*}
Hence, we can symmetrically prove that
\begin{align}\label{eq:Markov-lower-bound}
 \mathcal{N}(\bm{0},\bm{\Sigma} + \kappa \bm{I})\{\mathcal{A}\} -\mathbb{P}\left(\frac{\bm{S}_n}{\sqrt{n}} \in \mathcal{A}\right)
&\lesssim \sqrt{\kappa \log(nd)} \log_+d + \frac{\gamma^{\frac{3}{4}} + \beta^{\frac{3}{4}}}{\sqrt{\alpha}}\log_+^2d n^{-\frac{1}{4}}.
\end{align}
The theorem follows by combining \eqref{eq:Markov-Berry-Esseen-decompose-1}, \eqref{eq:Markov-Berry-Esseen-Gaussian-compare}, \eqref{eq:Markov-upper-bound} and \eqref{eq:Markov-lower-bound}.

\section{Proof of supportive lemmas and propositions}
\subsection{Proof of Inequality~\eqref{eq:prop.b1}}\label{app:proof-Stein-integral}
For simplicity, we denote
\begin{align*}
a = (1-e^{-2t})\varepsilon_1+e^{-2t}\varepsilon_k, \quad b=\varepsilon_k,
\end{align*}
and approach this integral through a change of variable
\begin{align*}
e^{-s} = \sqrt{\frac{a}{b}}\sin \theta.
\end{align*}
It can easily be verified that
\begin{align*}
e^{-s}\mathrm{d}s = -\sqrt{\frac{a}{b}}\cos \theta \mathrm{d}\theta, \quad \text{and} \quad a-be^{-2s} = a \cos^2 \theta.
\end{align*}
Therefore, the integral can be represented as 
\begin{align*}
\int_t^{\infty} \frac{e^{-3s}}{(a-be^{-2s})^{\frac{3}{2}}}\mathrm{d}s &= \int_t^{\infty} \frac{e^{-2s}}{(a-be^{-2s})^{\frac{3}{2}}} e^{-s}\mathrm{d}s \\ 
&= \int_0^{\theta^\star} \frac{\frac{a}{b}\sin^2 \theta}{(a \cos^2 \theta)^{\frac{3}{2}}} \cdot \sqrt{\frac{a}{b}}\cos \theta \mathrm{d}\theta \\ 
&= b^{-\frac{3}{2}} \int_0^{\theta^\star} \tan^2 \theta \mathrm{d}\theta \\ 
&= b^{-\frac{3}{2}} \cdot (\tan \theta - \theta)\bigg|_0^{\theta^\star} \\ 
&= b^{-\frac{3}{2}} \cdot (\tan \theta^\star - \theta^\star),
\end{align*}
where we let
\begin{align*}
e^{-t} = \sqrt{\frac{a}{b}}\sin \theta^\star, \quad \text{so} \quad \theta^\star = \arcsin\left(\sqrt{\frac{b}{a}}e^{-t}\right).
\end{align*}
Meanwhile, the difference between $\tan \theta^\star$ and $\theta^\star$ can be bounded by
\begin{align*}
\tan \theta^\star - \theta^\star &\leq \tan \theta^\star-\sin \theta^\star \\ 
&= \frac{\sqrt{\frac{b}{a}}e^{-t}}{\sqrt{1-\left(\sqrt{\frac{b}{a}}e^{-t}\right)^2}}-\sqrt{\frac{b}{a}}e^{-t} \\ 
&= \sqrt{\frac{e^{-2t}\varepsilon_k}{(1-e^{-2t})\varepsilon_1}} - \sqrt{\frac{e^{-2t}\varepsilon_k}{(1-e^{-2t})\varepsilon_1+e^{-2t}\varepsilon_k}} \\ 
&= \sqrt{\frac{e^{-2t}\varepsilon_k}{(1-e^{-2t})\varepsilon_1}} \left(1-\sqrt{\frac{(1-e^{-2t})\varepsilon_1}{(1-e^{-2t})\varepsilon_1+e^{-2t}\varepsilon_k}}\right) \\ 
& < \sqrt{\frac{e^{-2t}\varepsilon_k}{(1-e^{-2t})\varepsilon_1}} \left(1-\frac{(1-e^{-2t})\varepsilon_1}{(1-e^{-2t})\varepsilon_1+e^{-2t}\varepsilon_k}\right) \\ 
&= \sqrt{\frac{e^{-2t}\varepsilon_k}{(1-e^{-2t})\varepsilon_1}} \cdot \frac{e^{-2t}\varepsilon_k}{(1-e^{-2t})\varepsilon_1+e^{-2t}\varepsilon_k}.
\end{align*}
In combination, the integral is bounded by
\begin{align*}
&\int_t^\infty\frac{e^{-3s}}{((1-e^{2t})\varepsilon_1 + (e^{-2t}-e^{-2s})\varepsilon_k)^{\frac{3}{2}}} \mathrm{d}s \\ 
&= b^{-\frac{3}{2}} \cdot (\tan \theta^\star - \theta^\star) \\ 
&< \varepsilon_k^{-\frac{3}{2}} \cdot \sqrt{\frac{e^{-2t}\varepsilon_k}{(1-e^{-2t})\varepsilon_1}} \cdot \frac{e^{-2t}\varepsilon_k}{(1-e^{-2t})\varepsilon_1+e^{-2t}\varepsilon_k}\\ 
&= \frac{e^{-3t}}{\sqrt{(1-e^{-2t})\varepsilon_1}} \cdot \frac{1}{(1-e^{-2t})\varepsilon_1 + e^{-2t}\varepsilon_k}.
\end{align*}
\qed

\end{document}